\documentclass[11pt]{article}
\usepackage{amsfonts,mathrsfs,amssymb,amsthm,mathptm}
\usepackage{amsmath,amscd}
\usepackage{mathptm,pslatex}
\usepackage{color}
\usepackage[dvips]{graphicx}
\usepackage[all]{xy}
\usepackage{graphicx}
\usepackage{fancyhdr}
\usepackage{url}

\begin{document}
%\sloppy
\newtheorem{Def}{Definition}[section]
\newtheorem{Bsp}[Def]{Example}
\newtheorem{Prop}[Def]{Proposition}
\newtheorem{Theo}[Def]{Theorem}
\newtheorem{Lem}[Def]{Lemma}
\newtheorem{Koro}[Def]{Corollary}
\theoremstyle{definition}
\newtheorem{Rem}[Def]{Remark}

\newcommand{\add}{{\rm add}}
\newcommand{\con}{{\rm con}}
\newcommand{\gd}{{\rm gl.dim}}
\newcommand{\sd}{{\rm st.dim}}
\newcommand{\sr}{{\rm sr}}
\newcommand{\dm}{{\rm dom.dim}}
\newcommand{\cdm}{{\rm codomdim}}
\newcommand{\tdim}{{\rm dim}}
\newcommand{\E}{{\rm E}}
\newcommand{\Mor}{{\rm Morph}}
\newcommand{\End}{{\rm End}}
\newcommand{\ind}{{\rm ind}}
\newcommand{\rsd}{{\rm res.dim}}
\newcommand{\rd} {{\rm rd}}
\newcommand{\ol}{\overline}
\newcommand{\overpr}{$\hfill\square$}
\newcommand{\rad}{{\rm rad}}
\newcommand{\soc}{{\rm soc}}
\renewcommand{\top}{{\rm top}}
\newcommand{\pd}{{\rm pdim}}
\newcommand{\id}{{\rm idim}}
\newcommand{\fld}{{\rm fdim}}
\newcommand{\Fac}{{\rm Fac}}
\newcommand{\Gen}{{\rm Gen}}
\newcommand{\fd} {{\rm fin.dim}}
\newcommand{\Fd} {{\rm Fin.dim}}
\newcommand{\Pf}[1]{{\mathscr P}^{<\infty}(#1)}
\newcommand{\DTr}{{\rm DTr}}
\newcommand{\cpx}[1]{#1^{\bullet}}
\newcommand{\D}[1]{{\mathscr D}(#1)}
\newcommand{\Dz}[1]{{\mathscr D}^+(#1)}
\newcommand{\Df}[1]{{\mathscr D}^-(#1)}
\newcommand{\Db}[1]{{\mathscr D}^b(#1)}
\newcommand{\C}[1]{{\mathscr C}(#1)}
\newcommand{\Cz}[1]{{\mathscr C}^+(#1)}
\newcommand{\Cf}[1]{{\mathscr C}^-(#1)}
\newcommand{\Cb}[1]{{\mathscr C}^b(#1)}
\newcommand{\Dc}[1]{{\mathscr D}^c(#1)}
\newcommand{\K}[1]{{\mathscr K}(#1)}
\newcommand{\Kz}[1]{{\mathscr K}^+(#1)}
\newcommand{\Kf}[1]{{\mathscr  K}^-(#1)}
\newcommand{\Kb}[1]{{\mathscr K}^b(#1)}
\newcommand{\DF}[1]{{\mathscr D}_F(#1)}

\newcommand{\Kac}[1]{{\mathscr K}_{\rm ac}(#1)}
\newcommand{\Keac}[1]{{\mathscr K}_{\mbox{\rm e-ac}}(#1)}

\newcommand{\modcat}{\ensuremath{\mbox{{\rm -mod}}}}
\newcommand{\cmodcat}{\ensuremath{\mbox{{\rm -comod}}}}
\newcommand{\Modcat}{\ensuremath{\mbox{{\rm -Mod}}}}
\newcommand{\ires}{\ensuremath{\mbox{{\rm ires}}}}
\newcommand{\Stires}{\ensuremath{\mbox{{\rm Stires}}}}
\newcommand{\Stpres}{\ensuremath{\mbox{{\rm Stpres}}}}
\newcommand{\Spec}{{\rm Spec}}

\newcommand{\stmc}[1]{#1\mbox{{\rm -{\underline{mod}}}}}
\newcommand{\Stmc}[1]{#1\mbox{{\rm -{\underline{Mod}}}}}
\newcommand{\prj}[1]{#1\mbox{{\rm -proj}}}
\newcommand{\inj}[1]{#1\mbox{{\rm -inj}}}
\newcommand{\Prj}[1]{#1\mbox{{\rm -Proj}}}
\newcommand{\Inj}[1]{#1\mbox{{\rm -Inj}}}
\newcommand{\PI}[1]{#1\mbox{{\rm -Prinj}}}
\newcommand{\GP}[1]{#1\mbox{{\rm -GProj}}}
\newcommand{\GI}[1]{#1\mbox{{\rm -GInj}}}
\newcommand{\gp}[1]{#1\mbox{{\rm -Gproj}}}
\newcommand{\gi}[1]{#1\mbox{{\rm -Ginj}}}

\newcommand{\opp}{^{\rm op}}
\newcommand{\otimesL}{\otimes^{\rm\mathbb L}}
\newcommand{\rHom}{{\rm\mathbb R}{\rm Hom}\,}
\newcommand{\pdim}{\pd}
\newcommand{\Hom}{{\rm Hom}}
\newcommand{\Coker}{{\rm Coker}}
\newcommand{ \Ker  }{{\rm Ker}}
\newcommand{ \Cone }{{\rm Con}}
\newcommand{ \Img  }{{\rm Im}}
\newcommand{\Ext}{{\rm Ext}}
\newcommand{\StHom}{{\rm \underline{Hom}}}
\newcommand{\StEnd}{{\rm \underline{End}}}
\newcommand{\KK}{I\!\!K}
\newcommand{\gm}{{\rm _{\Gamma_M}}}
\newcommand{\gmr}{{\rm _{\Gamma_M^R}}}

\def\vez{\varepsilon}\def\bz{\bigoplus}  \def\sz {\oplus}
\def\epa{\xrightarrow} \def\inja{\hookrightarrow}

\newcommand{\lra}{\longrightarrow}
\newcommand{\llra}{\longleftarrow}
\newcommand{\lraf}[1]{\stackrel{#1}{\lra}}
\newcommand{\llaf}[1]{\stackrel{#1}{\llra}}
\newcommand{\ra}{\rightarrow}
\newcommand{\dk}{{\rm dim_{_{k}}}}

\newcommand{\holim}{{\rm Holim}}
\newcommand{\hocolim}{{\rm Hocolim}}
\newcommand{\colim}{{\rm colim\, }}
\newcommand{\limt}{{\rm lim\, }}
\newcommand{\Add}{{\rm Add }}
\newcommand{\Prod}{{\rm Prod }}
\newcommand{\pres}{\ensuremath{\mbox{{\rm pres}}}}
\newcommand{\app}{{\rm app }}
\newcommand{\Tor}{{\rm Tor}}
\newcommand{\Cogen}{{\rm Cogen}}
\newcommand{\Tria}{{\rm Tria}}
\newcommand{\Loc}{{\rm Loc}}
\newcommand{\Coloc}{{\rm Coloc}}
\newcommand{\tria}{{\rm tria}}
\newcommand{\Con}{{\rm Con}}
\newcommand{\Thick}{{\rm Thick}}
\newcommand{\thick}{{\rm thick}}
\newcommand{\Sum}{{\rm Sum}}

{\Large \bf
\begin{center}
Stable equivalences and homological dimensions
%Stable equivalences of centralizer matrix algebras and homological dimensions
%Stable equivalences and homological dimensions of centralizer matrix algebras.
\end{center}}

\medskip
\centerline{\textbf{Xiaogang Li} and \textbf{Changchang Xi}$^*$ }

\medskip
\centerline{Dedicated to the memory of Professor Roberto Mart\'inez-Villa (1942--2025)}
%\centerline{In memory of Professor Roberto Mart\'inez-Villa (19.10.1942-22.9.2025)}

\renewcommand{\thefootnote}{\alph{footnote}}
\setcounter{footnote}{-1} \footnote{ $^*$ Corresponding author.
Email: xicc@cnu.edu.cn; Fax: 0086 10 68903637.}
\renewcommand{\thefootnote}{\alph{footnote}}
\setcounter{footnote}{-1}
\renewcommand{\thefootnote}{\alph{footnote}}
\setcounter{footnote}{-1}

\begin{abstract}
As is known, every finite-dimensional algebra over a field is isomorphic to the centralizer algebra of \textbf{two} matrices. So it is fundamental to study first the centralizer algebra of a single matrix, called a centralizer matrix algebra. In this article, stable equivalences between centralizer matrix algebras over arbitrary fields are completely characterized in terms of a new type of equivalence relation on matrices. Moreover, stable equivalences of centralizer matrix algebras over any fields induce stable equivalences of Morita type, thus preserve dominant, finitistic and global dimensions. Our methods also show that the Alperin--Auslander/Auslander--Reiten conjecture holds true for stable equivalences between an arbitrary algebra and a centralizer matrix algebra over a common field.
\end{abstract}

{\footnotesize\tableofcontents\label{contents}}

\section{Introduction\label{Introduction}\label{sect1}}
In the representation theory of algebras and groups, stable equivalence is one of the fundamental categorical equivalences between different algebras. There is a large variety of literature on the subject. For example, the series of works of Auslander and Reiten on stable equivalence of dualizing $R$-varieties I-V (see \cite{ar5papers}, \cite{IR}). Stably equivalent algebras have many important invariants. One of them is that stably equivalent algebras have a one-to-one correspondence between non-projective, indecomposable modules (see \cite[Proposition 1.3, p.337]{ARS}), and another invariant is the dominant dimension of stably equivalent algebras without nodes and without semisimple summands (see \cite{MV2}). Despite great progresses made in the past decades, there are still many fundamental but  difficult questions on stable equivalences of algebras. For instance, how can we characterize stable equivalences of algebras (or special classes of algebras)? In this direction, there are only a few classes of algebras for which the problem is partially or completely solved. For example, Auslander and Reiten showed that an algebras with radical-square-zero is stably equivalent to hereditary algebras (see \cite{ar1975}, \cite{IR}). Also, stable equivalences between representation-finite self-injective algebras over an algebraically closed field induce derived equivalences (see \cite{asashiba2003}), and thus can be characterized by tensor products of bimodules or two-sided tilting complexes (see \cite{JR2}).

Another important question on stable equivalences is the not yet solved Auslander--Reiten conjecture \cite[Conjecture (5), p.409]{ARS} or Alperin--Auslander conjecture (\cite[Conjecture 2.5]{rouq}):

(ARC) Stably equivalent algebras over fields have the same number of non-isomorphic, non-projective simple modules.

Despite a lot of efforts made in the last decades, not much is known about this conjecture. For representation-finite algebras over algebraically closed fields, it was verified in \cite{MV1985}. This result was extended to Frobenius-finite algebras over algebraically closed fields under stable equivalences of Morita type in \cite{hx3}. Notably, the conjecture was reduced to self-injective algebras (see \cite{MV2}).

We aim to investigate stable equivalences and the related conjecture for an interesting class of algebras, called centralizer matrix algebras, that were studied initially by Georg Ferdinand Frobenius (see \cite{frob1877}), and have attracted considerable attentions in different areas (see \cite{cdfk, Datta, glr, stuart, xz2}).
Recently, centralizer matrix algebras have been investigated intensively in \cite{wz, xy, xz1, xz2}.
Now, let us recall their definition.

Let $n$ be a natural number, $R$ a field and $M_n(R)$ the full $n\times n$ matrix algebra over $R$ with the identity matrix $I_n$. For a nonempty subset $X$ of $M_n(R)$, we consider the centralizer algebra of $X$ in $M_n(R)$:
$$S_n(X,R):=\{a\in M_n(R)\mid ax=xa,\; \forall \; x\in X\}.$$
For simplicity, we write $S_n(c,R)$ for $S_n(\{c\},R)$. Clearly, $S_n(X,R)=\cap_{c\in C}S_n(c,R)$.

To understand $S_n(X,R)$ for a finite set $X$, it is enough to consider the case $|X|= 2$ by \cite[Lemma 1]{brenner1972}. Thus it is fundamental to understand first the centralizer algebra of a single matrix. Such an algebra $S_n(c,R)$ is called a \emph{centralizer matrix algebra} in this article.

Centralizer matrix algebras contribute an important part to general finite-dimensional algebras. This can be seen by a result, due to Sheila Brenner \cite{brenner1972, brenner1974}.

\smallskip
\textbf{Theorem.} (Brenner). \emph{Every finite-dimensional algebra over a field $R$ is isomorphic to the centralizer algebra of \textbf{two} matrices in a full  matrix algebra over $R$.}

\smallskip
Recently, a lot of new features of centralizer matrix algebras have been discovered. For instance, both (ARC) and the finitistic dimension conjecture are verified for centralizer matrix algebras over arbitrary fields (see \cite{xz3} and \cite{lx1}, respectively). Further, the Morita and derived equivalences of centralizer matrix algebras are completely described in  terms of $M$- and $D$-equivalences on matrices, respectively. Thus the categorical equivalences of centralizer matrix algebras are reduced to matrix equivalences in linear algebra \cite[Theorem 1.3]{lx1}.

The present article focuses on stable equivalences of centralizer matrix algebras.

\medskip
{\bf Question}: Given a field $R$ and two matrices $c\in M_n(R)$, $d\in M_m(R)$, how can we decide whether $S_n(c,R)$ and $S_m(d,R)$ are stably equivalent or not?

\smallskip
As is known, one of the main obstacles in dealing with stable equivalences is that there are not any criteria of stable equivalences in homological bimodules, while Morita and derived equivalences do have bimodule descriptions (see \cite{morita} and \cite{JR2}). Nevertheless, we will overcome this difficulty, introduce a new equivalence relation on matrices and provide surprisingly a complete characterization of stable equivalences between centralizer matrix algebras, as was done for Morita and derived equivalences in \cite{lx1}.

The new equivalence relation on square matrices is called an $S$-equivalence. It is defined purely in terms of linear algebra. We refer to Section \ref{sect3} for details.

Now, let us state our main results and their consequences more precisely.

\begin{Theo}\label{main1}
$(1)$ Let $R$ be a  field, $c$ and $d$ be square matrices over $R$ (may have different sizes). Then the centralizer matrix algebra of $c$ and the centralizer matrix algebra of $d$ are stably equivalent if and only if the matrices $c$ and $d$ are $S$-equivalent.

$(2)$ All algebras stably equivalent to a centralizer matrix algebra over a field have the same number of non-isomorphic, non-projective simple modules.
\end{Theo}

A polynomial of positive degree is said to be \emph{separable} if all of its irreducible factors have only simple roots in their splitting fields.
By Theorem \ref{main1}(1), the existence of a stable equivalence between centralizer matrix algebras can be decided by methods in linear algebra. Theorem \ref{main1}(2) generalizes slightly \cite[Theorem 1.4(1)]{xz3}, which says that (ARC) holds true for centralizer matrix algebras. Here, we do not require that both algebras considered are centralizer matrix algebras.

The strategy of the proof of Theorem \ref{main1} runs as follows: First, we show that  stable equivalences between centralizer matrix algebras preserve non-semisimple blocks. This is not true in general. Thus we have to consider the stable equivalence between blocks of centralizer matrix algebras. Second, we apply stable equivalences of Morita type to studying the stable equivalences of these blocks. This is carried out by investigating stable equivalences between the endomorphism algebras of generators over certain quotient algebras of the polynomial algebra $R[x]$. Finally, we lift a stable equivalence of Frobenius parts to a stable equivalence of algebras themselves.

As a consequence of Theorem \ref{main1}(1), we consider permutation matrices. Let $\Sigma_n$ denote the symmetric group of permutations on $\{1,2, \cdots, n\}$. For $\sigma\in \Sigma_n$, let $c_{\sigma}$ be the permutation matrix of $\sigma$. Given a prime number $p>0$,  we can define a $p$-regular part $r(\sigma)$ and a $p$-singular part $s(\sigma)$ of $\sigma$, and consider $r(\sigma)$ and $s(\sigma)$ as elements in $\Sigma_n$. As a convention, we set $r(\sigma)=\sigma$ when $p=0$. Let $c_{\sigma}$ denote the permutation matrix of $\sigma$ over $R$. We refer to Section \ref{stper} for more details.

Clearly, $S_n(c_{r(\sigma)},R)$ is a semisimple algebra (see \cite[Theorem 1.2(1)]{xz2}), and therefore its stable module category is trivial. But, for singular parts, we have the following.

\begin{Koro}\label{cor1.4}
Let $R$ be a field of characteristic $p\ge 0$, $\sigma\in \Sigma_n$ and $\tau\in \Sigma_m$. If $S_n(c_{\sigma},R)$ and $S_m(c_{\tau},R)$ are stably equivalent, then so are
$S_n(c_{s(\sigma)},R)$ and $S_m(c_{s(\tau)},R)$.
\end{Koro}

Finally, we point out an important property of centralizer matrix algebras over any fields.

\begin{Koro}\label{cor1.5}
$(1)$ Two centralizer matrix algebras over a common field are stably equivalent if and only if their non-semisimple parts are stably equivalent of Morita type.

$(2)$ Stably equivalent centralizer matrix algebras over a common field have the same global, finitistic and dominant dimensions.
\end{Koro}

\smallskip
The paper is outlined as follows. In Section \ref{sect2} we recall notions and terminologies. Also, we prepare some preliminaries for later proofs. In Section \ref{sect3} we introduce the $S$-equivalence relation on square matrices over a field.  In Section \ref{Pf} we prove all results mentioned in Section \ref{sect1}. During the course of proofs, we also generalize slightly a result in \cite[Theorem 1.4(1)]{xz3} about Auslander--Reiten conjecture on stable equivalences (see Proposition \ref{A-R}). Finally, we mention a couple of open problems related to results in this article. For example, how can we characterize stable equivalences between centralizer matrix algebras over a principal integral domain?

\section{Preliminaries}\label{sect2}

In this section we first fix notation and recall definitions and facts on both homological dimensions and stable equivalences, and then we prepare a few lemmas on modules over polynomial algebras.

\subsection{Definitions and notation\label{sect2.1}}
Throughout this paper, $R$ denotes a fixed field unless stated otherwise. Algebras always mean finite-dimensional unitary associative algebras over $R$, and modules mean  left modules.

Let $A$ be an algebra and $A^{\opp}$ the opposite algebra of $A$. By $A\modcat$ we denote the category of all finitely generated left $A$-modules, and by $A\prj$ (respectively, $A\inj$) the full subcategory of $A\modcat$ consisting of projective (respectively, injective) $A$-modules.
Let $A\modcat_{\mathscr{P}}$ stand for the full subcategory of $A\modcat$ consisting of modules without any nonzero projective direct summands.

For $M\in A\modcat$, let ${\add}(M)$ be the full subcategory of $A\modcat$ consisting of modules isomorphic to direct sums of finitely many indecomposable direct summands of $M$. Let $\ell(M)$ stand for the composition length of $M$. The basic module of $M$ is defined to be the direct sum of representatives of isomorphism classes of indecomposable direct summands of $M$, denoted by $\mathcal{B}(M)$. Let $M_{\mathscr{P}}$ denote the direct sum of all non-projective indecomposable summands of $M$. Thus $M/M\!_{\mathscr{P}}$ is projective. For $r\in \mathbb{N}$, we write $M^r$ for the direct sum of $r$ copies of $M$.

The module $_AM$ is called a \emph{generator} (or \emph{cogenerator}) for $A\modcat$ if $\add(_AA)\subseteq \add(M)$ (or $\add(D(A_A))$ $\subseteq \add(M)$, where $D: A\modcat\to A^{\opp}\modcat$ is the usual duality of $A$. By $-^*$ we denote the $A$-duality, that is $M^*:=\Hom_A(M,A)$.

Let $\pres(M)$ be the full subcategory of $A\modcat$ consisting of those modules $L$ such that there is an exact sequence:  $P_1\ra P_0\ra L\ra 0$ in $A\modcat$, with $P_0, P_1\in \add(M)$, and let $\app(M)$ be the full subcategory of $\pres(M)$ consisting of those modules $L$ that the induced map $\Hom_A(M,P_0)\to \Hom_A(M,L)$ is surjective.

For homomorphisms $f:X\to Y$ and $g: Y\to Z$ in $A\modcat$, we write $fg: X\to Z$ for their composition. This implies that the image of an element $x\in X$ under $f$ is denoted by $(x)f$. Thus $\Hom_A(X,Y)$ is naturally a left $\End_A(X)$- and right $\End_A(Y)$-bimodule.  The composition of functors between categories is written from right to left, that is, for two functors
$F:\mathcal{C}\ra \mathcal{D}$ and $G:\mathcal{D}\ra \Sigma$, we write $G\circ F$, or simply $GF$, for the composition of $F$ with $G$. The image of an object $X\in \mathcal{C}$ under $F$ is written as $F(X).$

Let $N\in A\modcat$. We write $\rad(M,N)$ for the set of homomorphisms $f\in \Hom_A(M,N)$ such that, for any homomorphisms $g:Z\ra M$ and $h:N\ra Z$, the composition $gfh: Z\to Z$ is not an automorphism of $Z$. Thus $\rad(M,M)$ is just the Jacobson radical of the endomorphism algebra $\End_A(M)$ of $M$.

The Nakayama functor $\nu_A:=D\Hom_A(-,A)\simeq D(A)\otimes_A-$ from $A\modcat$ to itself restricts to an equivalence from $\prj{A}$ to $A\inj$. An $A$-module $P$ is said to be $\nu$-\emph{stably projective} if $\nu^i_A{P}$ is projective for all $i\ge 0$. For example, if
$a^2=a\in A$ satisfies $\add(\nu_A(Aa))=\add(Aa)$, then $Aa$ is $\nu$-stably projective. In this case, the idempotent $a\in A$ is said to be \emph{$\nu$-stable}. Let $A$-stp denote the full subcategory of $A\modcat$ consisting of all $\nu$-stably projective $A$-modules. Clearly, there is a $\nu$-stable idempotent $e\in A$ such that $A$-stp = $\add(Ae)$. It is known that $eAe$ is a self-injective algebra (see \cite{hx2, MV1}).
Following \cite{hx2}, the algebra $eAe$ is called the \emph{Frobenius part} of $A$. This is uniquely determined by $A$ up to Morita equivalence.

If $A$ is a self-injective algebra and $M\in A\modcat$, then the Frobenius part of $\End_A(A\oplus M)$ is Morita equivalent to $A$.

For a class $\mathcal{C}$ of $A$-modules in $A$-mod, the number of modules in $\mathcal{C}$ always means the cardinality of the set of isomorphism classes of modules in $\mathcal{C}$.

Let $\eta:0\to X\to Y\to Z\to 0$ be an exact sequence in $A\modcat$. We say that $\eta$ has \emph{no split direct summands} if it has no split exact sequences as its direct summands.

\subsection{Homological dimensions and conjectures}
In this subsection, we recall the notions of dominant, finitistic and global dimensions.

Given an algebra $A$, the \emph{global dimension} of $A$, denoted $\gd(A)$, is defined by
$$\gd(A):=\mbox{sup}\{\mbox{proj.dim}(M)\mid M\in A\modcat\},$$
while the \emph{finitistic dimension} of $A$, denoted $\fd(A)$, is defined by

$$\fd(A):=\mbox{sup}\{\mbox{proj.dim}(M)\mid M\in A\modcat,\mbox{ proj.dim}(M)<\infty\}.$$

Clearly, $\fd(A)\le \gd(A)$. The well-known \emph{finitistic dimension conjecture} reads \cite{bass}:

\medskip
(FDC)  For any algebra $A$ over a field, $\fd(A)<\infty$.

\smallskip
This conjecture is a central problem in the representation theory of algebras. The validity of (FDC) implies the validity of several other homological conjectures (see \cite[Conjectures, p.409]{ARS}). Unfortunately, (FDC) is still open to date. In fact, only a few cases are verified. For example, it was verified for monomial algebras \cite{gk91} and algebras of radical cube-zero \cite{gz91}. In the last decades, there were several approaches to the conjecture. For instance, Igusa-Todorov's $\phi$- and $\psi$-dimensions \cite{ig05}, extensions of algebras \cite{xx}, removing arrows \cite{gps21} and delooping levels \cite{g22}.

Consider a minimal injective resolution of $_AA$:
$$0\lra {}_AA \lra I_0 \lra I_1 \lra \cdots \lra I_t \lra \cdots.$$ The \emph{dominant dimension} of $A$, denoted $\dm(A)$, is the maximal $t\in \mathbb{N}$ (or $\infty$) such  that $I_0, I_1,\cdots, I_{t-1}$ are projective. Related to dominant dimension, there is a well-known conjecture, called the \emph{Nakayama conjecture} (see \cite{n58}, or \cite[Conjectures, p.409]{ARS}):

(NC) If $\dm(A)=\infty,$ then $A$ is a self-injective algebra.

It is easy to see that the validity of (FDC) for an algebra $A$ implies the validity of (NC) for the algebra $A$. Note that both conjectures are open. However, we have verified them for centralizer matrix algebras \cite[Theorem 1.2]{lx1}.

\subsection{Stable equivalences and elimination of nodes\label{sect2.3}}
In this subsection we recall some basics on stable equivalences and a procedure of producing algebras without nodes.

By $A\stmc$ we denote the stable module category of $A$, which has the same objects as $A$-mod, but the morphism set $\underline{\Hom}_A(X,Y)$ of objects $X$ and $Y$ is the quotient of $\Hom_A(X,Y)$ modulo $\mathcal{P}(X,Y)$, the set of all homomorphisms from $X$ to $Y$ that factorize through  projective $A$-modules. If $M$ is a non-projective indecomposable $A$-module, then $\mathcal{P}(M,M)\subseteq \rad(\End_A(M))$, the Jacosbson radical of $\End_A(M)$.

Though we mainly work with finite-dimensional algebras over fields, we will state notions and facts in the most general form of Artin algebras, in order to be referred in other cases.

\begin{Def} \label{sta-equ}$(1)$ Two Artin algebras $\Lambda$ and $\Gamma$ over a commutative Artin ring $R$ are said to be \emph{stably equivalent} if there is an equivalence $F: \Lambda\stmc{} \ra \Gamma\stmc$ of $R$-linear categories.

$(2)$ A stable equivalence $F: \Lambda\stmc{} \ra \Gamma\stmc$ of Artin algebras
$\Lambda$ and $\Gamma$ \emph{preserves non-semisimple blocks} provided that, for non-projective indecomposable modules $M,N\in \Lambda\modcat_{\mathscr{P}}$, $F(M)$ and $F(N)$ lie in the same block of $\Gamma$ if and only if $M$ and $N$ lie in the same block of $\Lambda$.
\end{Def}

Let $F$ be a stable equivalence
between  $\Lambda$ and $\Gamma$. Then $F$ induces a one-to-one correspondence between $\Lambda\modcat_{\mathscr{P}}$ and  $\Gamma\modcat_{\mathscr{P}},$ but may \emph{not} preserve the numbers of non-semisimple blocks of algebras. This can be seen by the example in \cite[Example 4.7]{xz3}.

A special stable equivalence was introduced by Brou\'e for finite-dimensional algebras.

\begin{Def}{\rm \cite{MB}}\label{semt} Finite-dimensional algebras $A$ and $B$ over a field are \emph{stably equivalent of Morita type} if there exist bimodules $_AM_B$ and $_BN_A$ such that

$(1)$ $M$ and $N$ are projective as one-sided modules,

$(2)$ There are bimodule isomorphisms: $_AM\otimes_B N_A\simeq {}_AA_A\oplus P$ and $_BN\otimes_A M_B\simeq {}_BB_B\oplus Q$, where
$P$ is a projective $A$-$A$-bimodule and $Q$ is a projective $B$-$B$-bimodules.
\end{Def}

Clearly, the exact functor $N\otimes_A-: A\modcat{}\ra B\modcat$ induces a stable equivalence $N\otimes_A-:A\stmc{}\ra B\stmc$. By Definition \ref{semt}(1), proj.dim$(_BN\otimes_AX)\le $ proj.dim$(_AX)$ for all $X\in A\modcat$, where proj.dim$(X)$ denotes the projective dimension of $X$. Similarly, proj.dim$(_AM\otimes_BY)\le$ proj.dim$(_BY)$ for all $Y\in B\modcat$. Moreover, by  Definition \ref{semt}(2), we have proj.dim$(_AX)\le $ proj.dim$(X\oplus P\otimes_AX)$ = proj.dim$(M\otimes_BN\otimes_AX)\le$ proj.dim$(_BN\otimes_AX)\le $ proj.dim$(_AX)$. This implies that stable equivalences of Morita type preserve the global and finitistic dimensions.

An example of stable equivalences of Morita type is that derived equivalent self-injective algebras are stably equivalent of Morita type \cite{JR2}. Now we give another example of stable equivalences of Morita type.

\begin{Lem}\label{almm}
{\rm \cite{LX3}}
Let $B$ be a self-injective algebra and $X\in B\modcat$. Then the endomorphism algebras ${\End}_B(B\oplus X)$ and ${\End}_B(B\oplus \Omega_B(X))$ are stably equivalent of Morita type, where $\Omega_B(X)$ stands for the syzygy of $X$.
\end{Lem}

We discuss a special case of Lemma \ref{almm}: Let $B$ be a symmetric algebra. Consider the exact sequence of $B$-$B$-bimodules:
$$0\lra N\lra B\otimes_RB\lra B\lra 0.$$
Then $N$ is projective as a one-side $B$-module, and $N\otimes_B -$ induces a stable equivalence of Morita type between $B$ and itself. Claerly, $N\otimes_B- \simeq \Omega_B(-)$ on $B\stmc$.
% (see for instance \cite{As}).

Given a $B$-module $M$, we set $\Lambda:=\End_B(B\oplus M)$ and $\Gamma:=\End_B(B\oplus \Omega_B(M))$. Let $P:=\Hom_B(B\oplus M,B)$ and $Q:=\Hom_B(B\oplus \Omega_B(M), B)$. Then $_{\Lambda}P$ and $_{\Gamma}Q$ are projective modules. We may assume that $M$ is a basic module without nonzero projective summands. In \cite[Section3]{LX3}, a stable equivalence $F$ of Morita type is constructed between $\Lambda$ and $\Gamma$ such that the following diagram is commutative (up to natural isomorphism):
$$\xy
(0,15)*+{B\modcat}="a",
(40,15)*+{\pres(P)}="b",
(80,15)*+{\Lambda\modcat}="c",
(-20,8)*+{(\varepsilon)}="g",
(0,0)*+{B\modcat}="d",
(40,0)*+{\pres(Q)}="e",
(80,0)*+{\Gamma\modcat.}="f",
{\ar^{N\otimes_B-} "a";"d"},
{\ar^{F} "b";"e"},
{\ar^{F} "c";"f"},
{\ar^{P\otimes_B-}, "a";"b"},
{\ar^{}, "b";"c"},
{\ar^{Q\otimes_B-}, "d";"e"},
{\ar^{}, "e";"f"},
\endxy$$

\medskip
In the study of stable equivalences of Artin algebras, the notion of nodes plays an important role.

\begin{Def}{\rm \cite{MV1}} A non-projective, non-injective simple module over an Artin algebra is called a \emph{node} if the middle term of the almost split sequence starting at the simple module is projective.\end{Def}

By \cite[Lemma 1]{MV1}, a non-injective simple module $S$ of an Artin algebra is a node if and only if $S$ is not a composition factor of $\rad(Q)/\soc(Q)$ for any indecomposable projective module $Q$. Thus an Artin algebra has no nodes if and only if every non-projective, non-injective simple module is a composition factor of ${\rm rad}(P)/{\rm soc}(P)$ for some indecomposable projective module $P$.

Given an Artin algebra $\Lambda$, let $S$ be the direct sum of all non-isomorphic nodes of $\Lambda$, $I$ the trace of $S$ in $\Lambda$, and $J$ the left annihilator of $I$ in $\Lambda.$ Mart\'{i}nez-Villa showed in \cite[Theorem 2.10]{MV1} that an Artin algebra $\Lambda$ with nodes is stably equivalent to the triangular matrix algebra
$$\Lambda':=\left(
      \begin{array}{cc}    \Lambda/I & 0 \\     I & \Lambda/J \\   \end{array} \right)$$
without nodes. It is shown in \cite[Lemma 3.12(3)]{xz3} that $\Lambda$ and $\Lambda'$ have the same numbers of non-projective simple modules. We often say  that $\Lambda'$ is obtained from $\Lambda$ by eliminating nodes or $\Lambda'$ is the \emph{node-eliminated algebra} of $ \Lambda$.

\begin{Rem}\label{rmk2.11} Let $f(x)$ be an irreducible polynomial in $R[x]$, $A:=R[x]/(f(x)^2)$, $M(1):= R[x]/(f(x))$, $M(2):= A$ and $_AM := M(1)\oplus M(2)$. Clearly, $M(1)$ is the unique non-projective, indecomposable $A$-module. Let $C:= \End_A(M)$ be the Auslander algebra of $A$, and let $e_i\in C$ be the canonical projection of $M$ onto $M(i)$ for $i=1,2$. Then $C$ is a Nakayama algebra and has $2$ indecomposable projective modules  $P(1):=Ce_1$ and $P(2):=Ce_2$. The non-projective indecomposable $C$-modules are $S_1: = \top(P(1))$, $S_2:$ = top$(P(2))$ and the injective envelope $I(S_1)$ of $S_1$. Moreover, there are three almost split sequences of $C$-modules: $$0\lra S_1\lra I(S_1)\lra S_2\lra 0,\quad\quad0\lra S_2\lra P(1)\lra S_1\lra 0~\mbox{ and }$$
$$0\lra P(1)\lra P(2)\oplus S_1\lra I(S_1)\lra 0.$$Let $\Gamma^s_C$ be the quiver of $C$ obtained from the Auslander-Reiten quiver $\Gamma_C$ of $C$ by removing both projective vertices and all related arrows. Then $\Gamma^s_C$ is of the form $\bullet\ra\bullet\ra\bullet$.
By eliminating nodes, $A$ is stably equivalent to $A':=\begin{pmatrix}A/\rad(A) & 0 \\ \rad(A) &A/\rad(A)\\ \end{pmatrix}$,  and $C$ is stably equivalent to $C':=\begin{pmatrix} C/\soc(C) & 0 \\ \soc(C) & C/Ce_2C\\ \end{pmatrix}$. The algebra $A'$ has two indecomposable projective modules $\begin{pmatrix} A/\rad(A) \\ \rad(A)\\ \end{pmatrix}$ and $\begin{pmatrix} 0 \\ A/\rad(A)\\ \end{pmatrix}$. The first one is injective, while the second one is simple and isomorphic to a submodule of the first one. The algebra $C'$ has three indecomposable projective modules $\begin{pmatrix} \big(Ce_1+\soc(C)\big)/\soc(C)\\ \soc(Ce_1)\\ \end{pmatrix}$, $\begin{pmatrix} \big(Ce_2+\soc(C)\big)/\soc(C)\\ \soc(Ce_2)\\ \end{pmatrix}$ and $\begin{pmatrix} 0 \\ C/Ce_2C\\ \end{pmatrix}$. The second one is injective, while the third one is simple and isomorphic to a submodule of both the first and second ones. Note that the Frobenius parts of both $A'$ and $C'$ are $0$.
\end{Rem}

We recall the following results of Mart\'{i}nez-Villa in \cite[Proposition 1.5 and Theorems 1.7 and 2.6]{MV2}.

\begin{Lem}\label{exa}{\rm \cite{MV2}} Let $F:\Lambda\stmc$ $\ra \Gamma\stmc$ be a stable equivalence of Artin algebras $\Lambda$ and $\Gamma$ in which both algebras have neither nodes nor semisimple summands.

$(1)$ The functor $F$ provides a bijection $F':\mathscr{P}(\Lambda)_{\mathscr{I}}\ra \mathscr{P}(\Gamma)_{\mathscr{I}}$, which preserves simple projective modules in $\mathscr{P}(\Lambda)_{\mathscr{I}}$. The inverse of $F'$ is $H'$ induced from a quasi-inverse $H$ of the functor $F$.

$(2)$ The functor $F$ induces a stable equivalence between the Frobenius parts of $\Lambda$ and $\Gamma$.

$(3)$ Let $0\ra X\oplus Q_1\stackrel{f}\ra Y\oplus Q_2\oplus P\stackrel{g}\ra Z\ra 0$ be an exact sequence of $\Lambda$-modules without split exact sequences as its direct summands, where $X,Y,Z\in \Lambda\modcat_{\mathscr{P}}$, $Q_1,Q_2\in \mathscr{P}(\Lambda)_{\mathscr{I}}$ and $P$ is a projective-injective $\Lambda$-module. Then there is an exact sequence
$$\quad 0\lra F(X)\oplus F'(Q_1)\stackrel{f'}\lra F(Y)\oplus F'(Q_2)\oplus P'\stackrel{g'}\lra F(Z)\lra 0$$
in $\Gamma\modcat$, such that it has no split direct summands and $_{\Gamma}P'$ is projective-injective.
\end{Lem}

Here, $g'$ is a projective cover of $F(Z)$ if $g$ is a projective cover of $Z$, and $\mathscr{P}(\Lambda)_{\mathscr{I}}$ is the category of projective $\Lambda$-modules without nonzero injective direct summands.

As a consequence of Lemma \ref{exa}, we have the following result.
\begin{Lem}\label{exa1}
Let $F:\Lambda\stmc$ $\ra \Gamma\stmc$ be a stable equivalence of Artin algebras $\Lambda$ and $\Gamma$ in which both algebras have neither nodes nor semisimple summands, and let $P$ and $Q$ be the sums of projective-injective indecomposable $\Lambda$-modules and $\Gamma$-modules, respectively. Suppose $N\in \mathscr{P}(\Lambda)_{\mathscr{I}}$ and $M$ is a non-projective indecomposable module. If $M$ belongs to $\pres(P\oplus N)$, then $F(M)$ belongs to $\pres(Q\oplus F'(N))$.
\end{Lem}

{\it Proof.} Let $P(X)$ denote a projective cover of $X\in \Lambda\modcat$, and let $\eta(X)$ be the canonical exact sequence: $0\ra \Omega_\Lambda(X)\ra P(X)\ra X\ra 0$. Clearly, if $X\in \Lambda\modcat_{\mathscr{P}}$, then $\eta(X)$ has no split direct summands. Moreover, we have the following fact.

\smallskip
$(*)$ Let $Y\in \Lambda\prj$ and $X\in \Lambda\modcat$. Then $X\in \pres(Y)$ if and only if $P(X)$ and $ P(\Omega_\Lambda(X))$ lie in $\add(Y)$.

\smallskip
Now, let $M\in  \pres(P\oplus N)$ be non-projective and indecomposable. Then $P(M)$ and $P(\Omega_\Lambda(M))$ lie in $\add(P\oplus N)$. Since $\Omega_\Lambda(M)$ has no nonzero injective direct summands, we can write $\Omega_\Lambda(M)=(\bigoplus ^{s}_{i=1}X_i)\oplus Z$ with $X_i$ being non-projective indecomposable module for all $i\in [s]$ and $Z\in \mathscr{P}(\Lambda)_{\mathscr{I}}$. Therefore $P(\Omega_\Lambda(M))=P((\bigoplus ^{s}_{i=1}X_i)\oplus Z)=(\bigoplus^{s}_{i=1}P(X_i))\oplus Z\in \add(P\oplus N)$. This implies that $P(X_i)\in \add(P\oplus N)$ for all $i\in [s]$ and $Z\in \add(N)$. Since each of $\eta(X_1),\cdots, \eta(X_s)$ and $\eta(M)$ has no split direct summands, it follows from Lemma \ref{exa}(3) that there are exact sequences  $\eta(F(X_1)), \cdots, \eta(F(X_s))$ and $\eta(F(M))$ such that each of them has no split direct summands. From these exact sequences, we see that $P(F(X_i))\in \add(Q\oplus F'(N))$ for all $i\in [s]$, $P(F(M))\in \add(Q\oplus F'(N))$ and $\Omega_\Gamma(F(M))=F(\bigoplus ^{s}_{i=1}X_i)\oplus F'(Z)=(\bigoplus^s_{i=1}F(X_i))\oplus F'(Z)$. Clearly, $F'(Z)\in \add(Q\oplus F'(N))$ since $Z\in \add(N)$. Thus $P(\Omega_\Gamma(F(M)))=(\bigoplus^s_{i=1}P(F(X_i)))\oplus F'(Z)\in \add(Q\oplus F'(N))$. Hence $F(M)\in \pres(Q\oplus F'(N))$ by $(*)$. $\square$

\subsection{Modules over quotients of polynomial algebras}
In this subsection we show some facts on stably equivalent endomorphism algebras of generators over quotients of $R[x]$.

Throughout this subsection, let $[n]$ be the set $\{1,2,\cdots, n\}$, $f(x)$ a fixed irreducible polynomial in $R[x]$ and $A:=R[x]/(f(x)^n)$.

Clearly, $A$ is a local, symmetric Nakayama algebra.
Moreover, $A$ has $n$ indecomposable modules $M(i):=R[x]/(f(x)^i)$, $i\in [n].$ For convenience, we set $M(0)=0$. For $n\ge 2$, let $\Gamma_{n-1}:=\{M(i)\mid i\in [n-1]\}\subseteq A\modcat_{\mathscr{P}}$ be the set of representatives of isomorphism classes of non-projective indecomposable $A$-modules.

Suppose that $F:\stmc{A}  \to A\stmc$ is a stable equivalence. Then $F$ induces an action $\ol{F}$ on $\Gamma_{n-1}$, namely, for  $M\in \Gamma_{n-1}$, $\ol{F}(M)$ is the unique module in $\Gamma_{n-1}$ such that $\ol{F}(M)\simeq F(M)$ in $A\stmc$. Clearly, for the stable equivalence $\Omega_A: \stmc{A}  \to A\stmc$, we have $\overline{\Omega}_A(M(i))= M(n-i)$, where $\Omega_A$ is the syzygy operator of $A$.

\begin{Lem}{\rm \cite[Lemma 2.10]{lx1}}\label{self} Let $n\ge 2$. If $F$ is a stable equivalence between $A$ and itself, then the induced action $\ol{F}$ on $\Gamma_{n-1}$ coincides with either the identity action or $\overline{\Omega}_A$.
\end{Lem}

For a positive integer $m$ and an irreducible polynomial $g(x)\in R[x]$, if $A\simeq R[x]/(g(x)^m)$ as $R$-algebras, then $n=LL(R[x]/(f(x)^n))=LL(R[x]/(g(x)^m))=m$. Moreover, the indecomposable $R[x]/(g(x)^n)$-module $R[x]/(g(x)^t)$ is isomorphic to the $A$-module $M(t)$ for $t\in [n]$.

Let $R$ be a commutative local ring, and let $R_i:=R/\rad^i(R)$ for all $i\in \mathbb{N}$. For $j\leq i$, there is the natural surjective homomorphism $\pi_{ij}: R_i\ra R_j$ of rings.  The inverse limit of the inverse system $(R_i, \pi_{ij})$ is denoted by $\hat{R}$. The commutative local ring $R$ is said to be {\it complete} if the rings $R$ and $\hat{R}$ are isomorphic.

We mention the following structure theorem of commutative complete local ring by Cohen \cite[Theorem 9]{Cohen}.
\begin{Lem}\label{Co}
Let $R$ be a commutative complete local ring, and suppose that $R$ has the same characteristic as $R/\rad(R)$. Then $R$ contains a unitary subring $S$ such that the restriction of the natural map $R\ra R/\rad(R)$ to $S$ is an isomorphism of rings.
\end{Lem}
As consequences  of  Lemma \ref{Co}, we have the following two facts:

$(1)$ Any commutative Artin local rings are always complete because their maximal ideals are nilpotent ideals;

$(2)$ If the ring $R$ in Lemma \ref{Co} is an algebra over a field, then the subring $S$ can be chosen to be a subalgebra over the same field. This follows from the proof of \cite[Theorem 9]{Cohen}.

Thus it follows from the above facts (1)-(2) and Lemma \ref{Co} that $A$ always contains a subalgebra $K$ which is isomorphic to $R[x]/(f(x))$. This implies that we can exclude the restriction that $R$ is perfect in \cite[Lemma 2.14, Corollary 2.15]{lx1}. So the improved version of the two results now read as follows.

\begin{Lem}\label{poi}
$(1)$ Let $g(x)\in R[x]$ be an irreducible polynomial and $m>0$. If $A\simeq R[x]/(g(x)^m)$ as $R$-algebras, then $m=n$. Moreover, for $t\in [n]$, the indecomposable $R[x]/(g(x)^n)$-module $R[x]/(g(x)^t)$ is isomorphic to the $A$-module $M(t)$.

$(2)$ $A$ can be viewed as a $K$-algebra and $A\simeq K[x]/(x^n)$ as algebras over $K:=R[x]/(f(x))$.

$(3)$ If $g(x)\in R[x]$ is irreducible such that $A$ is stably equivalent to $R[x]/(g(x)^m)$ for $m\ge 2$, then $A\simeq R[x]/(g(x)^m)$ as $R$-algebras.
\end{Lem}

Now we consider the endomorphism algebras of generators over $A$.

\begin{Lem}\label{node} Let $M$ be a generator for $A\modcat$, that is $_AA\in \add(M)$, and $\Lambda:=\End_A(M)$.

$(1)$ If $n\geq 2$, then every simple $\Lambda$-module is neither projective nor injective.

$(2)$ $\Lambda$ has nodes if and only if $n = 2$.

$(3)$ The Frobenius part of $\Lambda$ is Morita equivalent to $A$.
\end{Lem}

{\it Proof.} Since the Nakayama algebra $A$ is local and symmetric, any indecomposable direct summand of $M$ is isomorphic to a submodule of $_AA$. Let $F$ be the functor $\Hom_A(M,-): \add(_AM)\to \prj{\Lambda}$ and $P:=F(_AA).$ Clearly, $F$ is an equivalence and $P$ is an indecomposable projective-injective $\Lambda$-module with $\soc(P)\simeq \top(P)$ since $A$ is a local symmetric Nakayama algebra. The left exactness of $F$ implies that any indecomposable projective $\Lambda$-module is isomorphic to a submodule of $P$.

(1) Suppose $n\geq 2$. Then  the indecomposable, projective-injective $\Lambda$-module $P$ is not simple. To show that every indecomposable projective $\Lambda$-module is not simple, we show that $\Hom_{\Lambda}(X,Y)\ne 0$ for any indecomposable projective $\Lambda$-modules $X$ and $Y$. In fact, let $X'$ and $Y'$ be indecomposable direct summands of $_AM$ such that $F(X')=X$ and $F(Y')=Y$. Since $A$ is a local algebra, we  always have $\Hom_A(X',Y')\ne 0$ (for example, a nonzero homomorphism from the top of $X'$ to the socle of $Y'$). Thus $$\Hom_\Lambda(X,Y)=\Hom_{\Lambda}(F(X'),F(Y'))\simeq \Hom_A(X',Y')\ne 0.$$
This implies that no simple $\Lambda$-module is projective. Otherwise, $\Lambda$ would have only one simple projective module, this would contradict to the fact that $P$ is not simple. Since the Nakayama functor $\nu:\prj{\Lambda}\to \Lambda\inj$ is an equivalence, we deduce that, for any indecomposable injective $\Lambda$-modules $U$ and $V$,
$$\Hom_\Lambda(U,V)\simeq \Hom_\Lambda(\nu^{-1}(U),\nu^{-1}(V))\ne 0.$$ This shows that no simple $\Lambda$-modules are injective.

(2) If $n=1$, then the algebra $A$ is simple, and therefore the algebra $\Lambda$ is semisimple. Thus $\Lambda$ has no nodes by definition.

If $n=2$, then $\Lambda$ is Morita equivalent to either $A$ or the Auslander algebra $C$ of $A$. Note that $C$ is a Nakayama algebra with $2$ indecomposable projective $C$-modules $P_1$ and $P_2=P$ of lengths $2$ and $3$, respectively. There is an almost split sequence $0\to \top(P_2)\to P_1\to \top(P_1)\to 0$, which shows that $\Lambda$ has a node.

Finally, we consider $n\geq 3$. Let $S$ be a non-projective, non-injective simple $\Lambda$-module, we show that $S$ is a composition factor of $\rad(P)/\soc(P).$ Actually, let $Q$ be the projective cover of $S$ with $Q=\Hom_A(M,Y)$ for an indecomposable $_AY\in \add(M)$. Then $Q$ is a non-simple submodule of $P$ with $\soc(Q)=\soc(P)$. If $Q\neq P$, then $Q\subset \rad(P)$ and $0\ne Q/\soc(Q)$ is a submodule of $\rad(P)/\soc(P)$. Thus, as a composition factor of $Q/\soc(Q)$, $S$ is a composition factor of $\rad(P)/\soc(P)$. If $Q=P$, then the multiplicity of $\top(P)$ in $P$ is at least $LL(A)\ge 3$, and therefore $S=\top(P)$ is a composition factor of $\rad(P)/\soc(P)$.  Hence $\Lambda$ has no nodes.

(3) This is clear.
$\square$

\begin{Lem} \label{sta-down} Suppose that $M$ is a generator for $A\modcat$. Let $B:= R[x]/(g(x)^m)$ for some irreducible polynomial $g(x)$ and $m>0$, and let $N$ be a generator for $B\modcat$. If $\End_A(M)$ and $\End_B(N)$ are stably equivalent, then so are $A$ and $B$.
\end{Lem}

{\it Proof.} Let $\Lambda:=\End_A(M)$ and $\Gamma:=\End_B(N)$. By Lemma \ref{node}(3), the Frobenius part of $\Gamma$ is Morita equivalent to $B$.
Suppose that $F: \Lambda\stmc{}\ra \Gamma\stmc$ is an equivalence of $R$-categories. Then $\Lambda$ is semisimple if and only if so is $\Gamma$. This implies that $m=1$ if and only $n=1$. Thus we may assume both $n\geq 2$ and $m\geq 2$. Under this assumption, we are led to considering the two cases.

$(1)$ $m\ge 3$ and $n\ge 3$. Then both $\Lambda$ and $\Gamma$ have no nodes by Lemma \ref{node}(2). Thus it follows from Lemma \ref{exa}(2) that $A$ and $B$ are stably equivalent.

$(2)$ $m=2$ or  $n=2$. In this case, we claim $m=n=2$. Suppose contrarily $m\neq n$. Then one of $m$ and $n$ is at least $3$. We may assume $n=2$ and $m\geq 3$. According to Lemma \ref{node}(2), $\Lambda$ has nodes, but $\Gamma$ has no nodes. We consider the node-eliminated algebra $\Lambda'$ of $\Lambda$. Clearly, $\Lambda'$ is stably equivalent to $\Lambda$ and has the Frobenius part equal to $0$ (see Remark \ref{rmk2.11}). Note that both $\Lambda'$ and $\Gamma$ have no semisimple summands. Thus, by Lemma \ref{exa}(2), the Frobenius part of $\Lambda'$ and the Frobenius part $B$ of $\Gamma$ are stably equivalent, and therefore  $\stmc{B}=0$. However, for $m\geq 2$, the $B$-module $R[x]/(g(x))$ is non-projective, and therefore $\stmc{B}\neq 0$, a contradiction. This shows $m=n=2$.

The Auslander algebra of $A:=R[x]/(f(x)^2)$ (respectively, $B:=R[x]/(g(x)^2)$) is a Nakayama algebra with $3$ non-projective indecomposable modules, one of which is injective of length $2$ and the other two are simple. Since stably equivalent algebras have the same number of non-projective indecomposable modules, it follows that

(i) either $\Lambda$ is Morita equivalent to $A$, and  $\Gamma$ is Morita equivalent to $B$; or

(ii) $\Lambda$ is Morita equivalent  to the Auslander algebras of $A$, and  $\Gamma$ is Morita equivalent to the Auslander algebras of $B$.

To complete the proof, we only need to deal with (ii). By using Hom-functor, we see that the endomorphism algebra of any simple $\Lambda$-module (respectively, $\Gamma$-module) is isomorphic to $R[x]/(f(x))$ (respectively, $R[x]/(g(x))$). By the definition of stable endomorphism algebras of modules, we see that the stable endomorphism algebra of any simple $\Lambda$-module (respectively, $\Gamma$-module) is isomorphic to $R[x]/(f(x))$ (respectively, $R[x]/(g(x))$). Clearly, there exists a simple $\Lambda$-module, say $T$, such that $F(T)$ is simple. Since $F$ is an equivalence of categories, it follows that $R[x]/(f(x))\simeq\StEnd_\Lambda(T)\simeq \StEnd_\Gamma(F(T))\simeq R[x]/(g(x))$ as algebras. This yields that $A$ and $B$ are stably equivalent since $R[x]/(f(x))$ and $R[x]/(g(x))$ are the only non-projective indecomposable $A$- and $B$-module, respectively. $\square$

\smallskip
Finally, we point out the following facts (see \cite[Lemma 2.19]{lx1} and \cite[Theorem 1.2(1)]{xz2}).

\begin{Lem}\label{iso-pr}
$(1)$ For $c\in M_n(R)$, there are isomorphisms of $R$-algebras: $$S_n(c,R)\simeq S_n(c^{tr},R)\simeq S_n(c,R)^{\opp}\simeq \End_{R[c]}(R^n),$$where $c^{tr}$ denotes the transpose of the matrix $c$.

$(2)$ If $R$ is of characteristic $p\geq 0$ and $\sigma\in \Sigma_n$, then $S_n(c_{\sigma},R)$ is semisimple if and only if $\sigma=r(\sigma)$, where $c_{\sigma}$ is the permutation matrix of $\sigma$ over $R$ and $r(\sigma)$ is the $p$-regular part of $\sigma$.
%\end{Lem}
\end{Lem}

\section{New equivalence relation on square matrices\label{sect3}}
In this section we introduce an $S$-equivalence on all square matrices over a field.

Let $R[x]$ the polynomial algebra over a field $R$ in one variable $x$.
For polynomials $f(x)$ and  $g(x)$ of positive degree, we write $f(x)\mid g(x)$ or $f(x)\le g(x)$ if $f(x)$ divides $g(x)$, that is, $g(x)=f(x)h(x)$ with $h(x)\in R[x]$.  Then the set of monic polynomials in $R[x]$ of positive degree is actually a partially ordered set with respect to this polynomial divisibility.

The elementary divisors of $c\in M_n(R)$ is defined to the elementary divisors of the matrix $xI_n-c\in M_n(R[x])$, and all elementary divisors of $c$ form a multiset. Let $\mathcal{E}_c$ be the set of elementary divisors of $c$ (thus $\mathcal{E}_c$ has no duplicate elements), and let $\mathcal{R}_c$ be the subset of $\mathcal{E}_c$ consisting of all reducible elementary divisors of $c$ that are maximal with respect to $\le$.

For $f(x)\in \mathcal{R}_c$, we define the set $P_c(f(x))$ of power indices of $f(x)$ by
$$P_c(f(x)):=\{i\ge 1 \mid  p(x)^i\in \mathcal{E}_c \mbox{ with } p(x) \mbox{ an irreducible factor of } f(x)\}.$$

Let $\mathbb{Z}_{>0}$ be the set of all positive integers and $s\in \mathbb{Z}_{>0}$. Given a set $T:=\{ m_1, m_2, \cdots, m_s\}\subset\mathbb{Z}_{>0}$ with $m_1>m_2>\cdots >m_s$, we associate it with a set $\mathcal{J}_T=:\{m_1,m_1-m_2,\cdots,m_1-m_s\}$. For a finite subset $H\subset\mathbb{Z}_{>0}$, $H=\mathcal{J}_T$ if and only if $T=\mathcal{J}_H$.

Now, we introduce a new equivalence relation on square matrices, which is different from the ones in \cite[Definition 1.1]{lx1}.
\begin{Def}\label{s-equ}
Matrices $c\in M_n(R)$ and $d\in M_m(R)$ are $S$-equivalent, denoted $c\stackrel{S}\sim d$, if there is a bijection $\pi: \mathcal{R}_c\ra \mathcal{R}_d$ such that, for $f(x)\in \mathcal{R}_c$, the two conditions hold:

$(1)$ $R[x]/(f(x))\simeq R[x]/((f(x))\pi)$ as algebras over $R$, and

 $(2)$ either $P_c(f(x))= P_d((f(x))\pi)$ or $P_c(f(x))=\mathcal{J}_{P_d((f(x))\pi)}$.
\end{Def}

This is an equivalence relation on all square matrices over $R$.
Let us give examples to illustrate $S$-equivalent matrices.

\begin{Bsp}{\rm
 Let $J_n(\lambda)$ denote the $n\times n$ Jordan block with the eigenvalue $\lambda\in R$.

$(1)$ We take $c=J_3(0)\oplus J_1(0)\oplus J_1(1)\in M_5(R)$ and $d=J_3(1)\oplus J_2(1)\in M_5(R)$, where $\oplus$ means taking block diagonal matrix. Then $\mathcal{E}_c=\{x^3,x,x-1\}$ and $\mathcal{E}_d=\{(x-1)^3, (x-1)^2\}$. By definition, $\mathcal{R}_c=\{x^3\}$ and  $\mathcal{R}_d=\{(x-1)^3\}$. Further, $P_c(x^3)=\{1,3\}$, $P_d((x-1)^3)=\{2,3\}$ and $\mathcal{J}_{P_d((x-1)^3)}=\{1,3\}$. Let $\pi:\mathcal{R}_c\ra \mathcal{R}_d$ be the map, $x^3\mapsto (x-1)^3$. Then $c\stackrel{S}\sim d$ by Definition \ref{s-equ}. This also shows that $c_1=J_3(0)\oplus J_1(0)\in M_4(R)$ is $S$-equivalent to $d$.

(2) The $S$-equivalence is different from the $D$-equivalence introduced in \cite{lx1}. On the one hand, the matrices $c$ and $d$ in (1) are $S$-equivalent, but not $D$-equivalent by \cite[Definition 1.1]{lx1}.
On the other hand, let
$c:=J_5(0)\oplus J_4(0)\oplus J_1(0)\in M_{10}(R)$ and $d:=J_5(0)\oplus J_2(0)\oplus J_1(0)\in M_8(R)$. By \cite[Example 4.4]{lx1}, $c$ and $d$ are $D$-equivalent, but not $S$-equivalent.
}
\end{Bsp}

Clearly, the similarity of matrices preserves $S$-equivalence, that is, if $c$ is similar to $c_1$ and if $d$ is similar to $d_1$, then $c\stackrel{S}\sim d$ implies $c_1\stackrel{S}\sim d_1$.

\section{Stable equivalences of centralizer matrix algebras\label{Pf}}
This section is devoted to proofs of main results and their corollaries by following the strategy mentioned in the introduction.

For $c\in M_n(R)$,  let $m_c(x)$ be the minimal polynomial of $c$ over $R$ and $A_c:=R[x]/(m_c(x))$. We write
$$m_c(x):=\prod^{l_c}_{i=1} f_i(x)^{n_i} \mbox{  for } n_i\ge 1 \; \mbox{ and } \;U_i:=R[x]/(f_i(x)^{n_i})\mbox{  for } 1\le i\le  l_c$$ where all $f_i(x)$ are distinct irreducible (monic) polynomials in $R[x]$. It is known that $U_i$ is a local, symmetric Nakayama $R$-algebra for $1\le i\le l_c$, and
$$A_c\simeq U_1\times U_2\times\cdots\times U_{l_c}.$$

Since $A_c\simeq R[c]$ and the $R[c]$-module $R^n=\{(a_1,a_2,\cdots,a_n)^{tr}\mid a_i\in R, 1\le i\le n\}$ can be regarded as an $A_c$-module, we can decompose the $A_c$-module $R^n$, according to the blocks of $A_c$, in the following way:
$$(\star)\quad \quad R^n\simeq \bigoplus^{l_c}_{i=1}\bigoplus^{s_i}_{j=1} \; R[x]/(f_i(x)^{e_{ij}})$$
as $R[x]$-modules, where $e_{ij}$ are positive integers. Note that $\{\{f_i(x)^{e_{ij}}\mid i\in [l_c],  j\in [s_i]\}\}$ is the multiset of all elementary divisors of $c$ (see \cite[Chapter 4]{Co}, where $(\star)$ is stated in terms of invariant subspaces of a linear transformation).

Let $M_i:= \bigoplus^{s_i}_{j=1} \; R[x]/(f_i(x)^{e_{ij}})$ be the sum of indecomposable direct summands of $R^n$  belonging to the block $U_i$, and $A_i:=\End_{U_i}(M_i)$ for $i\in  [l_c]$. Then all algebras $A_i$ are indecomposable.

Clearly, for $i\in [l_c]$, we have $P_c(f_i(x)^{n_i})=\{e_{ij}\mid j\in [s_i]\}$ and $\mathcal{B}(M_i)\simeq \bigoplus_{r\in {P_c(f_i(x)^{n_i})}} R[x]/(f_i(x)^r)$ as $U_i$-modules. Since $R^n$ is a faithful $M_n(R)$-module, $R^n$ is also a faithful $R[c]$-module, and therefore $R^n$ is a generator for $R[c]\modcat$, and $M_i$ is a faithful $U_i$-module for $i\in [l_c]$.
$$S_n(c,R)\simeq \End_{R[c]}(R^n)\simeq \End_{A_c}(\bigoplus_{i=1}^{l_c}M_i)=\prod^{l_c}_{i=1}{\End}_{U_i}(M_i)= \prod^{l_c}_{i=1}A_i .$$
This is a decomposition of blocks of $S_n(c,R)$. Moreover, a block $A_i$ is semisimple if and only if $n_i=1$. As the $R[c]$-module $R^n$ is a generator, we see that the bimodule $_{R[c]}R^n_{S_n(c,R)}$ has the double centralizer property, that is, $\End_{S_n(c,R)^{\opp}}(R^n_{_{S_n(c,R)}})=R[c]$.

The following lemma is an immediately consequence of $(\star)$.

\begin{Lem} \label{bijection}
There is a bijection $\pi$ from $\mathcal{E}_c$ to the set of non-isomorphic, indecomposable direct summands of the $A_c$-module $R^n$, sending $h(x)$ to the $A_c$-module $R[x]/(h(x))$ for $h(x)\in \mathcal{E}_c$.
\end{Lem}

Similarly, for $d\in M_m(R)$, we write $m_d(x)=\prod^{l_d}_{j=1} g_j(x)^{m_j} \; \mbox{  for } m_j\ge 1,$ where all $g_j(x)$ are distinct irreducible (monic) polynomials in $R[x]$. Let $A_d := R[x]/(m_d(x))$ = $V_1\times \cdots \times V_{l_d}$ with $V_j:=R[x]/(g_j(x)^{m_j})$ for $j\in [l_d]$. Clearly, $V_j$ is a local, symmetric Nakayama $R$-algebra. By the canonical isomorphism $A_d\simeq R[d]$ of algebras, we have a decomposition of the $A_d$-module $R^m$:
$$(\star\star)\quad \quad R^m\simeq \bigoplus^{l_d}_{i=1}\bigoplus^{t_i}_{j=1} \; R[x]/(g_i(x)^{f_{ij}}),$$
where the multiset $\{\{g_i(x)^{f_{ij}}\mid i\in [l_d], j\in [t_i]\}\}$ is the multiset of elementary divisors of $d$. This is also a decomposition of $R[x]$-modules. Let $N_i:= \bigoplus^{t_i}_{j=1} \; R[x]/(g_i(x)^{f_{ij}})$ be the sum of indecomposable direct summands of $R^n$  belonging to the block $V_i$, and $B_i:=\End_{V_i}(N_i)$ for $i\in [l_d]$. Then all algebras $B_i$ are indecomposable. Moreover, a block $B_i$ is semisimple if and only if $m_i=1$.

For $i\in [l_d]$, we have $P_d(g_i(x)^{m_i})=\{f_{ij}\mid j\in [t_i]\}$ and $\mathcal{B}(N_i)\simeq \bigoplus_{r\in {P_d(g_i(x)^{m_i})}} R[x]/(g_i(x)^r)$ as $V_i$-modules, and the following isomorphisms hold:

$$S_m(d,R)\simeq \End_{R[d]}(R^m)\simeq \End_{A_d}(\bigoplus_{i=1}^{l_d}N_i)=\prod^{l_d}_{i=1}{\End}_{V_i}(N_i)= \prod^{l_d}_{i=1}B_i .$$
This is a decomposition of blocks of $S_m(d,R)$ and the bimodule $_{R[d]}R^m_{S_m(d,R)}$ has the double centralizer property, that is, $\End_{S_m(d,R)^{\opp}}(R^m_{_{S_m(d,R)}})=R[d]$.

\subsection{Stable equivalences between blocks of centralizer matrix algebras}

In this subsection we study behaviour of blocks of centralizer matrix algebras under stable equivalences.

Let $\Lambda$ be an Artin algebra. We denote by $\mathscr{P}(\Lambda)_{\mathscr{I}}$  the set of representatives of isomorphism classes of projective $\Lambda$-modules without nonzero injective direct summands.

In general, stably equivalent algebras may have different numbers of non-semisimple blocks (see \cite[Example 4.7]{xz3}). In the following, we will show that stable equivalences between centralizer matrix algebras over a field do preserve non-semisimple blocks (see Definition \ref{sta-equ}).

First, we fix some notation. We have a block decomposition $S_n(c,R)\simeq\End_{R[c]}(R^n)\simeq \prod^{l_c}_{i=1}{\End}_{U_i}(M_i)=\prod^{l_c}_{i=1}A_i$. Let $A$ be the direct sum of all non-semisimple blocks of $S_n(c,R)$. Now, we partition the blocks $A_i$ of $A$ in the following way.

(1) $\mathcal{S}_{A,\geq 2}:=\{A_1, A_2, \cdots, A_{a_1}\}$ consists of the blocks $A_i$ with $n_i\ge 3$ and having at least 2 non-injective indecomposable projective modules for $1\le i\le a_1$.

(2) $\mathcal{S}_{A,1}:=\{A_{a_1+1}, A_{a_1+2}, \cdots, A_{a_2}\}$ consists of the blocks $A_i$ with $n_i\ge 3$ and having only 1 non-injective indecomposable projective module  for $a_1< i\le a_2$.

(3) $\mathcal{S}_{A,0}:=\{A_{a_2+1}, A_{a_2+2}, \cdots, A_{a_3}\}$ consists of the blocks $A_i$ with $n_i\ge 3$ and all indecomposable projective $A_i$-modules being injective for $a_2<i\le a_3$.

(4) $\mathcal{S}_{A}:=\{A_{a_3+1}, A_{a_3+2},\cdots, A_{a_4}\}$ consists of the blocks $A_i$ with $n_i=2$ for $a_3<i\le a_4$, where $0\le a_1\le a_2\le a_3\le a_4\le l_c.$

\smallskip
Note that $A=\prod_{i=1}^{a_4}A_i$ and $n_i$ is the Loewy length of both $U_i$ and the center of $A_i$. For an Artin algebra $\Lambda$, we denote by $\Lambda'$ the node-eliminated algebra of $\Lambda$. Thus $\Lambda$ and $\Lambda'$ are stably equivalent, and the latter has no nodes.

\medskip
Let $\widetilde{A}:=\prod^{a_3}_{i=1} A_i\; \times \prod_{a_3< i\le a_4} (A_i)'.$ Then we have the following result.

\begin{Lem}\label{replace} The algebra $\widetilde{A}$ has neither nodes nor semisimple direct summands. Moreover, there exists a stable equivalence
$F_A: A\stmc{}\ra \widetilde{A}\stmc$ such that $F_A$ preserves non-semisimple blocks of the algebras and its restriction to a block $A_i$ is the identity functor for $i\in [a_3]$.
\end{Lem}

{\it Proof.} Recall that $U_i= R[x]/(f_i(x)^{n_i})$ and $A_i = \End_{U_i}(M_i)$ for $i\in [l_c]$. For $i\in [a_3],$ it follows from Lemma \ref{node}(1)-(2) that $A_i$ has neither nodes nor projective simple modules. For $a_3<j\le a_4$, the node-eliminated algebra $(A_j)'$ of $A_j$ has no nodes and is stably equivalent to $A_j$. Clearly, $\widetilde{A}$ does not have any semisimple direct summands. Now it is easy to get a desired stable equivalence $F_A$ between $A$ and $\widetilde{A}$. $\square$

\medskip
For $S_m(d,R)$, let $B$ be the sum of its non-semisimple blocks. Similarly, we have a partition of blocks of $B$: there are natural numbers $0\le b_1\le b_2\le b_3\le b_4\le l_d$, such that the blocks of $B$ are partitioned as $\{B_1, \cdots, B_{b_1}\}\cup \{ B_{b_1+1},\cdots, B_{b_2}\} \cup \{B_{b_2+1},\cdots,B_{b_3}\}\cup \{B_{b_3+1},\cdots, B_{b_4}\}$. This partition of blocks of $B$ has the corresponding properties (1)-(4) as the partition of blocks of $A$.

\medskip
Let $\widetilde{B}:=\prod^{b_3}_{j=1} B_j \; \times \prod_{b_3< j\le b_4} (B_j)'$. Then, by Lemma \ref{replace}, $\widetilde{B}$
has neither nodes nor semisimple direct summands, and there is a stable equivalence $F_B: B\stmc{}\ra \widetilde{B}\stmc$ such that $F_B$ preserves non-semisimple blocks and its restriction to a block $B_j$ is the identity functor for  $j\in [b_3]$.

\medskip
From now on, we assume that \textbf{there is a stable equivalence $F$} between $S_n(c,R)$ and $S_m(d,R)$.

The functor $F$ restricts to a stable equivalence between $A$ and $B$. Thus $H:=F_B\circ F\circ F_A^{-1}:\widetilde{A}\stmc$ $\to$ $\widetilde{B}\stmc$ is a stable equivalence. Let $J: \widetilde{B}\stmc$ $\to \widetilde{A}\stmc$ be a quasi-inverse of $H$. As in Lemma \ref{exa}(1), $H$ and $J$ induce bijections $H': \mathscr{P}(\widetilde{A})_{\mathscr{I}}\ra \mathscr{P}(\widetilde{B})_{\mathscr{I}}$ and $J': \mathscr{P}(\widetilde{B})_{\mathscr{I}}\ra \mathscr{P}(\widetilde{A})_{\mathscr{I}}$, respectively, such that their compositions are identity maps.

\medskip
Let $\mathcal{S}_{A,\geq 1}:=\mathcal{S}_{A,\geq 2}\cup \mathcal{S}_{A,1}$ and $\mathcal{S}_{B,\geq 1}:=\mathcal{S}_{B,\geq 2}\cup \mathcal{S}_{B,1}$. In the following, we will show that $H$ restricts to a stable equivalence between blocks in $\mathcal{S}_{A,\geq 1}$ and blocks in $\mathcal{S}_{B,\geq 1}$.

\begin{Lem}\label{bij}
The bijection $H'$ induces a bijection between  $\mathcal{S}_{A,\geq 1}$ and  $\mathcal{S}_{B,\geq 1}$.
\end{Lem}

{\it Proof.} Recall the following facts:

(i) $U_i$ is a symmetric, local Nakayama algebra, and therefore any indecomposable $U_i$-module is isomorphic to a submodule of $U_i$ for $i\in [l_c]$.

(ii) $M_i$ is a generator for $U_i\modcat$. Let $A_i:=\End_{U_i}(M_i)$. Then any indecomposable projective $A_i$-module is isomorphic to a submodule of $P_i:=\Hom_{U_i}(M_i,U_i)$ for $i\in [l_c]$. This follows from (i) and the left exactness of Hom-functors.

(iii) $P_i$ is projective-injective by the isomorphism $\nu_{A_i}\Hom_{U_i}(M_i, U_i) \simeq \Hom_{U_i}(M_i,\nu_{U_i}U_i)$ (see \cite[Remark 2.9 (2)]{hx3}). In particular, $P_i$ is the unique (up to isomorphism) indecomposable projective-injective $A_i$-module for $i\in [l_c]$.

Now, let $i\in [a_2]$. The proof will be proceeded by the following two steps.

\smallskip
Step 1: We show that $H'$ induces a bijection from  $\mathcal{S}_{A,\geq 2}$ to $\mathcal{S}_{B,\geq 2}$. Thus $a_1=b_1$.

\smallskip
In fact, let $A_i\in \mathcal{S}_{A,\geq 2}$, that is, $i\in [a_1], \; n_i\ge 3$ and $A_i$ has at least $2$ non-isomorphic, non-injective indecomposable projective modules $P_{i1}$ and $P_{i2}$ which are not simple by Lemma \ref{node}(1). Then there exist $2$ indecomposable direct summands $M_{i1}$ and $M_{i2}$ of the $U_i$-module $M_i$ such that $P_{ik} = \Hom_{U_i}(M_i,M_{ik})$ for $1\le k\le 2$. As $U_i$ is a local Nakayama algebra, we may assume that $M_{i1}$ is isomorphic to a proper submodule of $M_{i2}$. It then follows from the left exactness of the Hom-functor $\Hom_{U_i}(M_i,-)$ that $P_{i1}$ is isomorphic to a submodule of $P_{i2}$. Hence there is
an exact sequence of $A_i$-modules
$$0\lra P_{i1}\lra P_{i2}\lra P_{i2}/P_{i1}\lra 0$$
with $P_{i2}/P_{i1}$ indecomposable. By Lemma \ref{exa}(3), there is an exact sequence of $\widetilde{B}$-modules
$$(*)\quad 0\lra H'(P_{i1})\lra H'(P_{i2})\oplus P'_1\lra H(P_{i2}/P_{i1})\lra 0,$$
with $P'_1$ being a projective-injective $\widetilde{B}$-module.
We show that both $H'(P_{i1})$ and $H'(P_{i2})$ lie in the same block of $\widetilde{B}$.

Suppose contrarily that $H'(P_{i1})$ and $H'(P_{i2})$ lie in different blocks of $\widetilde{B}$. Then $\Hom_{\widetilde{B}}(H'(P_{i1}), H'(P_{i2})) = 0$. It then follows from $(*)$ that the non-projective indecomposable module $H(P_{i2}/P_{i1})$  contains an isomorphic copy of $H'(P_{i2})$ as a direct summand. This is a contradiction. Therefore $H'(P_{i1})$ and $H'(P_{i2})$ lie in the same block of $\widetilde{B}$.

Now, suppose that $H'(P_{i1})$ and $H'(P_{i2})$ lie in a common block $C$ of $\widetilde{B}$. By Remark \ref{rmk2.11}, a block $(B_j)'$ of $\widetilde{B}$, with $b_3<j\le b_4$, always has at most 2 non-isomorphic, non-injective indecomposable projective modules, one of which is simple. It then follows from Lemma \ref{exa}(1) that $C\not\simeq (B_j)'$ as algebras for $b_3<j\le b_4$. Hence $C$ is a block of the form $B_j$ for some $j\in [b_3]$. Since the two non-injective indecomposable projective modules $H'(P_{i1})$ and $H'(P_{i2})$ lie in $C$, we deduce that $C$ is a block of the form $B_j$ for some $j\in [b_1]$ by our partition of blocks. Hence $H'$ sends all non-injective, indecomposable projective $A_i$-modules to the ones belonging to the block $B_j$ with $j\in [b_1]$. Similarly, by considering the the quasi-inverse $J$ of $H$ and the bijection map $J'$ which is the inverse of $H'$, we know that the non-injective indecomposable projective $B_j$-modules are mapped by $J'$ into modules belonging to the block $A_i$. Thus $H'$ restricts to a one-to-one correspondence between $\mathscr{P}(A_i)_{\mathscr{I}}$ and $\mathscr{P}(B_j)_{\mathscr{I}}$. This implies that $H'$ induces a bijection from $\mathcal{S}_{A,\geq 2}$ to $\mathcal{S}_{B,\geq 2}$, such that the corresponding blocks have the same number of non-isomorphic, non-injective indecomposable projective modules. Thus $a_1=b_1$.

\smallskip
Step 2:  We prove that $H'$ induces a bijection between $\mathcal{S}_{A,1}$ and $\mathcal{S}_{B,1}$. Thus $a_2=b_2$.

\smallskip
Actually, let $A_i\in \mathcal{S}_{A,1}$, that is, $a_1< i\le a_2$ and $A_i$ has only $1$ non-injective indecomposable projective module, say $Q_i$. If $H'(Q_i)$ lies in a block $(B_j)'$ for $b_3<j\le b_4$, then, by Remark \ref{rmk2.11}, $H'(Q_i)$ has a simple projective submodule, say $\widetilde{P_j}$. With a similar argument as in Step 1, we deduce that $Q_i$ and the simple projective module $J'(\widetilde{P_j})$ lie in the same block $A_i$. Note that $Q_i$ is not simple and $Q_i\not\simeq J'(\widetilde{P_j})$. Thus the block $A_i$ contains at least $2$ non-injective indecomposable projective modules. This is a contradiction. Thus it follows from Step 1 that $H'(Q_i)$ belongs to a block $B_j\in \mathcal{S}_{B,1}$. Conversely, for $B_j\in \mathcal{S}_{B,1}$, there is a unique non-injective indecomposable projective $B_j$-module, say $R_j$. As $J'$ is the inverse of $H'$, we see that $J'(R_j)$ belongs to a block $A_i\in \mathcal{S}_{A,1}$. So $H'$ induces a bijection between the set of blocks in $\mathcal{S}_{A,1}$ and  the set of blocks in $\mathcal{S}_{B,1}$. Thus $a_2=b_2.$ $\square$

Having established a bijection between  $\mathcal{S}_{A,\geq 1}$ and $\mathcal{S}_{B,\geq 1}$, we now show that the corresponding blocks are stably equivalent.

\begin{Lem}\label{st-eq-blocks} Let $H': \mathcal{S}_{A,\geq 1}\ra \mathcal{S}_{B,\geq 1}$ be the bijection in Lemma \ref{bij}. If $H'(A_i)=B_j$, then the functor $H$ restricts to a stable equivalence between $A_i$ and $B_j$.
\end{Lem}

{\it Proof.} For convenience, we define an operator $\hat{H}$ from the set of representatives of isomorphism classes of indecomposable $A_i$-modules to that of indecomposable $B_j$-modules:
$$\hat{H}(X)=\begin{cases}H'(X) & \mbox{ if } X \mbox{ is projective,}\\ H(X) & \mbox{ otherwise.}\end{cases}
$$

By Lemma \ref{bij}, we may assume $H'(A_i)=B_i$ for $i\in [a_2]$. To show that $H$ restricts to a stable equivalence between $A_i$ and $B_i$, we only need to show that, for any non-projective indecomposable $A_i$-module $M$, $H(M)$ belongs to the block $B_i$, and that, for any non-projective indecomposable $B_i$-module $N$, $J(N)$ belongs to the block $A_i$.

Indeed, let $i\in [a_2]$ and $M\in A_i\modcat$ be a non-projective indecomposable module. Then $H(M)$ is indecomposable. Now, we prove that $H(M)$ lies in the block $B_i$. We divide the proof into two parts (a) and (b) below.

(a) We show that, for the unique indecomposable projective-injective $A_i$-module $P_i$, $\hat{H}(\rad(P_i))$ lies in the block $B_i$.

As a submodule of the indecomposable injective module $P_i$, $\rad(P_i)$ has simple socle and therefore is indecomposable. Recall that each non-injective indecomposable projective $A_i$-module is isomorphic to a submodule of $P_i$. This implies that each non-injective, indecomposable projective $A_i$-module is isomorphic to a submodule of $\rad(P_i)$. By our partitions of blocks, $A_i$ has at least one  non-injective, indecomposable projective $A_i$-module, say $P$. In particular, $P$ is isomorphic to a submodule of $\rad(P_i)$. If $\rad(P_i)$ is projective, then it follows from $H'(A_i)=B_i$ that $\hat{H}(\rad(P_i))=H'(\rad(P_i))$ lies in the block $B_i$, and therefore (a) follows.

Now assume that $\rad(P_i)$ is not projective. Then $\rad(P_i)$ and $P$ are non-isomorphic, and so there exists an exact sequence
$$\quad 0\lra P\stackrel{\iota}\lra \rad(P_i)\stackrel{\eta}\lra \rad(P_i)/P\lra 0.$$
Since $\rad(P_i)$ is indecomposable, $\rad(P_i)/P$ does not have any nonzero projective direct summands.
Applying Lemma \ref{exa}(3) to this sequence, we get an exact sequence of $\widetilde{B}$-modules:
$$0\lra H'(P)\stackrel{\iota'}\lra H(\rad(P_i))\oplus P'_2\stackrel{\eta'}\lra H(\rad(P_i)/P)\lra 0,$$
with $P'_2$ a projective-injective $\widetilde{B}$-module. Due to $H'(A_i)=B_i$, we see that $H'(P)$ belongs to the block $B_i$.

Contrarily, suppose that $H(\rad(P_i))$ does not belong to the block $B_i.$ Then $\Hom_{\widetilde{B}}(H'(P), H(\rad(P_i)))=0$, and therefore $\Img(\iota')\subset P'_2.$ Thus $H(\rad(P_i)/P)\simeq (H(\rad(P_i))\oplus P'_2)/\Img(\iota')\simeq H(\rad(P_i))\oplus P'_2/\Img(\iota')$. Applying the quasi-inverse $J$ of $H$ to this isomorphism, we deduce that $\rad(P_i)$ is isomorphic to a direct summand of $\rad(P_i)/P$. Thus $P=0$. This contradicts to $P$ being indecomposable.
%that $\ell(\rad(P_i))\leq \ell(\rad(P_i)/P)<\ell(\rad(P_i))$, a contradiction.
Thus $\hat{H}(\rad(P_i))=H(\rad(P_i))$ belongs to the block $B_i.$

(b) We show that $H(M)$ belongs to the block $B_i$.

Indeed, let $P(M)$ be a projective cover of the indecomposable $A_i$-module $M$. Then there is a canonical exact sequence of $A_i$-modules
$ 0\ra \Omega(M)\stackrel{f}\ra P(M)\stackrel{g}\ra M\ra 0,$
which has no split exact sequences as its direct summands.
Write $\Omega(M)=L_1\oplus L_2$, with $L_1\in {A_i\modcat}_{\mathscr{P}}$ and $L_2\in \mathscr{P}(A_i)_{\mathscr{I}}$, and $P(M)= Q\oplus (P^t_i)$, with $Q\in \mathscr{P}(A_i)_{\mathscr{I}}$ and $t\in \mathbb{N}.$ By Lemma \ref{exa}(3), we get an exact sequence of $\widetilde{B}$-modules
$$(\ddag)\quad 0\lra H(L_1)\oplus H'(L_2)\stackrel{f'}\lra H'(Q)\oplus P'_3\stackrel{g'}\lra H(M)\lra 0,$$
such that $P'_3$ is a projective-injective $\widetilde{B}$-module and ($\ddag$) has no split direct summands. We have to consider the two cases: $Q\neq 0$ and $Q=0$.

(1) Assume $Q\neq 0$. We claim that $H(M)$ lies in the block $B_i$. Suppose contrarily that $H(M)$ does not belong to $B_i$. Then it follows from $H'(A_i)=B_i$ that $H'(Q)$ and $H(M)$ lie in different blocks, and therefore $\Hom_{\widetilde{B}}(H'(Q),H(M))=0$. Thus $ \Img(f')=H'(Q)\oplus (\Img(f')\cap P'_3)$. Since $H'(Q)$ is projective, there exists a homomorphism $h':H'(Q)\ra H(L_1)\oplus H'(L_2)$ such that $h'f'p_1=id_{H'(Q)}$, where $p_1$ is the projection from $H'(Q)\oplus P'_3$ to $H'(Q)$. This implies that the exact sequence
$$0\lra \Img(h')\stackrel{(f'|_{\Img(h')})p_1}\lra H'(Q)\lra 0\lra 0$$is a split direct summand of ($\ddag$). This is a contradiction, and therefore $H(M)$ lies in the block $B_i$.

(2) Assume $Q=0$. Then $H'(Q)=0$. According to the blocks of $\widetilde{B}$, we write $P'_3=\bigoplus^l_{k=1} X_k$ such that $X_i$ and $X_j$ are in different blocks for $i\ne j$, but all indecomposable direct summands of $X_i$ belong to a common block of $\widetilde{B}$. We shall show that all indecomposable direct summands of $P'_3$ lie in the same block of $\widetilde{B}$, that is, $l=1$.
Contrarily, suppose $l\geq 2$. By $(\ddag)$, we have $\Img(f')=\bigoplus^l_{k=1} Y_k$, where $Y_k$ is a submodules of $X_k$ for $k\in [l]$. Since $H(L_1)\oplus H'(L_2)$ contains no projective-injective direct summands, each $Y_k$ is a proper submodules of $X_k$. Therefore $H(M)\simeq P'_3/\Img(f')\simeq \bigoplus^l_{k=1} \; X_k/Y_k$ is decomposable. This contradicts to the fact that $H(M)$ is indecomposable. Thus $l=1$ and all indecomposable direct summands of $P'_3$,  including $P'_3$ iteself, lie in a common block of $\widetilde{B}$.

Suppose that $\Omega(M)$ has a direct summand isomorphic to $\rad(P_i)$. Then $f'$ in $(\ddag)$ restricts to an injective homomorphism from  $\hat{H}(\rad(P_i))$ to $P'_3$. Therefore $\hat{H}(\rad(P_i))$ and $P'_3$ lie in the same block $B_i$ by (a). Further, $H'(Q)=0$ implies $\Hom_{\widetilde{B}}(P'_3,H(M))\neq 0$ in $(\ddag)$. This yields that $H(M)$ and $P'_3$ lie in the same block $B_i$.

Suppose that $\Omega(M)$ has no direct summands isomorphic to $\rad(P_i)$. Under the assumption $Q=0$, we have $P(M)=P_i^{t}$. We consider the following exact sequence of $A_i$-modules
$$(\sharp)\quad 0\lra \Omega(M)\stackrel{f}\lra (\rad(P_i))^{t}\stackrel{g\mid_{(\rad(P_i))^{t}}}\lra \rad(M)\lra 0.$$
Clearly,
the module $\rad(M)$ is not projective, that is, $(\rad(M))_{\mathscr{P}}\neq 0$. By the assumption on $\Omega(M)$, any non-zero split direct summand of ($\sharp$) is of the form $$0\lra 0\lra (\rad(P_i))^{k}\stackrel{g\mid_{(\rad(P_i))^{k}}}\lra Y\lra 0$$ for a positive integer $k\leq t$ and a direct summand $Y$ of $\rad(M)$. Deleting split direct summands of the exact sequence ($\sharp$), we obtain an exact sequence of $A_i$-modules
$$0\lra \Omega(M)\stackrel{f_0}\lra (\rad(P_i))^{r}\stackrel{g_0}\lra X\lra 0,$$where $r\leq t$ and $X$ is a nonzero direct summand of $ (\rad(M))_{\mathscr{P}}$. Recall that $\Omega(M)=L_1\oplus L_2$ for $L_1\in {A_i\modcat}_{\mathscr{P}}$ and $L_2\in \mathscr{P}(A_i)_{\mathscr{I}}$.
According to Lemma \ref{exa}(3), we get an exact sequence of $\widetilde{B}$-modules:
$$0\lra H(L_1)\oplus H'(L_2)\stackrel{f'_0}\lra \big(\hat{H}(\rad(P_i))\big)^{r}\oplus P'_4\stackrel{g'_0}\lra H(X)\lra 0,$$where $P'_4$ is a projective-injective $\widetilde{B}$-module.

Suppose that $\Omega(M)$ contains no indecomposable direct summand $L$ such that $\hat{H}(L)$ lies in the block $B_i$. Then it follows from (a) that $\Hom_{\widetilde{B}}\big( H(L_1)\oplus H'(L_2),(\hat{H}(\rad(P_i)))^{r}\big)=0$ and so the image of $f'_0$ lies in $P'_4$. This implies that $H(X)$ contains $\big(\hat{H}(\rad(P_i))\big)^{r}$ as a direct summand. Since $H(X)$ does not have any non-zero projective direct summands and $H'(A_i)=B_i$, we deduce that $\rad(P_i)$ is not  projective and $\big(\hat{H}(\rad(P_i))\big)^{r}=\big(H(\rad(P_i))\big)^{r}$. Using the quasi-inverse $J$ of $H$, we deduce that $X$ contains $\rad(P_i)^{r}$ as a direct summand. This, together with the surjective homomorphism $g_0$, implies that $X\simeq \rad(P_i)^{r}$ and  $\Omega(M)=0$. Therefore $M\simeq P(M)$ is projective. This contradicts to our choice of $M$. Hence $\Omega(M)$ has an indecomposable direct summand $L$ such that $\hat{H}(L)$ lies in the block $B_i$. Then it follows from the non-split exact sequence ($\ddag$) that the modules $H(M)$, $P'_3$ and $\hat{H}(L)$ lie in the same block $B_i$.
Thus we have proved that $H(M)$ lies in the block $B_i$.

Similarly, for $i\in [b_2]$ and a non-projective, indecomposable $B_i$-module $N$, we see that $J(N)$ belongs to the block $A_i$ with $i\in [a_2]$. Hence $H$ induces a stable equivalence between $A_i$ and $B_i$ for $i\in [a_2]$. $\square$

\medskip
Next, we show that $H$ restricts to a stable equivalence between blocks in $\mathcal{S}_{A,0}\cup \mathcal{S}_{A}$ and blocks in $\mathcal{S}_{B,0}\cup \mathcal{S}_{B}$.

For this purpose, we recall the Auslander--Reiten quiver of stable module category of an Artin algebra. Given an Artin algebra $\Lambda$, let $\Gamma_{\Lambda}$ denote the Auslander--Reiten quiver of $\Lambda$, and let $\Gamma^s_{\Lambda}$ be the subquiver  of $\Gamma_{\Lambda}$ obtained by removing all projective vertices from $\Gamma_{\Lambda}$. For a local, symmetric Nakayama algebra $\Lambda_0:= R[x]/(f(x)^n)$ with $f(x)$ an irreducible polynomial, the quiver $\Gamma^s_{\Lambda_0}$ of $\Lambda_0$ is a connected quiver with $n-1$ vertices and $\tau(M(i))\simeq M(i)$ for $i\in [n-1]$. Here $\tau$ is the Auslander--Reiten translation. Moreover, there is an arrow $M(i)\to M(j)$ in $\Gamma^s_{\Lambda_0}$ if and only if there is an arrow from $M(j)\to M(i)$ in $\Gamma^s_{\Lambda_0}$ for $i,j\in [n-1]$.

From Lemma \ref{st-eq-blocks}, we see that $H$ restricts to a stable equivalence between $\prod^{a_2}_{i=1} A_i$ and $\prod^{b_2}_{j=1} B_j$, which preserves non-semisimple blocks, and therefore $H$ restricts to a stable equivalence between
$$C:=\prod_{a_2<i\leq a_3} A_i\; \times \prod_{a_3< i\le a_4} (A_i)' \; \mbox{ and } \; D:=\prod_{b_2<j\leq b_3} B_j\; \times \prod_{b_3< j\le b_4} (B_j)',$$where $(A_i)'$ and $(B_j)'$ are the node-eliminated blocks of $A_i$ and $B_j$, respectively. These products are actually decompositions of blocks of $C$ and $D$, respectively

In the following, we shall show that the restriction of $H$ between $C$ and $D$ preserves non-semisimple blocks as well as node-eliminated blocks, that is, for non-projective indecomposable $C$-modules $M$ and $N$, both $H(M)$ and $H(N)$ lie in a common node-eliminated block of $D$ if and only if both $M$ and $N$ lie in a common node-eliminated block of $C$.

First, we note the following immediate consequence of \cite[Lemma 1.2(d), p.336]{ARS}.

\begin{Lem}\label{preserve} Suppose that $G$ is a stable equivalence between algebras $\prod^{s}_{i=1}C_i$ and $\prod^{t}_{j=1}D_j$, where $C_i$ and $D_j$ are indecomposable non-semisimple algebras for $i\in [s]$ and $j\in [t]$. Suppose that the quivers $\Gamma^s_{C_i}$ and $\Gamma^s_{D_j}$ are connected for all $i\in [s]$ and all $j\in [t]$. Then $G$ preserves non-semisimple blocks, and therefore $s=t$.
\end{Lem}

\begin{Lem}\label{no-p-i} The functor $H$ induces a stable equivalence between $C$ and $D$, preserving non-semisimple blocks and node-eliminated blocks.
\end{Lem}

{\it Proof.}  We have seen that $H$ restricts to a stable equivalence between $C$ and $D$. For convenience, we denote this restriction by $H_C$.

For a block $A_i$ of $C$ with $a_2<i\leq a_3$, we have $A_i\in \mathcal{S}_{A,0}$, that is, $n_i\geq 3$ and all indecomposable projective $A_i$-modules are injective. Thus the $U_i$-module $M_i$ is projective and $A_i= \End_{U_i}(M_i)$ is Morita equivalent to $U_i=R[x]/(f_i(x)^{n_i})$. In particular, $A_i$ is a symmetric Nakayama algebra with $n_i-1$ non-projective, indecomposable modules. Thus, for $a_2<i\leq a_3$, $\Gamma^s_{A_i}$ is a connected quiver with two arrows between any two connected vertices.
For a block $(A_i)'$ of $C$ with $a_3<i\leq a_4$, we have $A_i\in \mathcal{S}_{A}$. In this case, $\Gamma^s_{(A_i)'}$ is either $\bullet$ or of the form $\bullet\ra\bullet\ra\bullet$ by Remark \ref{rmk2.11}.

Similarly, for a block $B_j$ of $D$ with $b_2<j\leq b_3$, the quiver $\Gamma^s_{B_j}$ is a connected quiver with two arrows between any two connected vertices; and  for $b_2<j\leq b_3$, the quiver $\Gamma^s_{(B_j)'}$ is either $\bullet$ or of the form $\bullet\ra\bullet\ra\bullet$.

According to Lemma \ref{preserve}, we know that $H_C$ preserves non-semisimple blocks. To complete the proof, it suffices to show that

(i) a block $A_i$ with $a_2< i \le a_3$ is never stably equivalent to a block $(B_j)'$ with $b_3< j \le b_4$, and

(ii) a block $(A_i)'$ with $a_3<i\leq a_4$ is never stably equivalent to a block $B_j$ with $b_2<j\leq b_3$.

Actually, $\Gamma^s_{A_i}$ (respectively, $\Gamma^s_{B_j}$) contains a subquiver of the form $\xymatrix{
\bullet\ar@<2.5pt>[r]^{} &\bullet \ar@<2.5pt>[l]^{}}$ for $a_2<i\leq a_3$ (respectively, $b_2<j\leq b_3$), while $\Gamma^s_{(A_i)'}$ (respectively, $\Gamma^s_{(B_j)'}$) contains no such subquivers for $a_3<i\leq a_4$ (respectively, $b_3<j\leq b_4$). According to \cite[Lemma 1.2(d), p. 336]{ARS}, we infer that (i) and (ii) hold true.  $\square$

\medskip
Combining the above discussions, we arrive at the following result.
\begin{Prop}\label{sta}
Suppose that there is a stable equivalence $F$ between $S_n(c,R)$ and $S_m(d,R)$. Then $F$ preserves non-semisimple blocks. Moreover,  if $F$ induces a stable equivalence between non-semisimple blocks $A_i$ and $B_j$ with $i,j\in [a_4]$, then $n_i=m_j$.
\end{Prop}

{\it Proof.} Suppose that there is a stable equivalence $F$ between $S_n(c,R)$ and $S_m(d,R)$.
Note that $A$ and $B$ denote the sums of non-semisimple blocks of $S_n(c,R)$ and $S_m(d,R)$, respectively, and  that $F_A: A\stmc{}\ra \widetilde{A}\stmc$ and $F_B: B\stmc{}\ra \widetilde{B}\stmc$ are stable equivalences constructed by eliminating nodes (see Lemma \ref{replace}).

By Lemmas \ref{bij} and \ref{no-p-i}, the stable equivalence $H=F_B\circ F\circ F_A^{-1}:\widetilde{A}\stmc{}\ra \widetilde{B}\stmc$ preserves non-semisimple blocks. Since both $F_A$ and $F_B$ preserve non-semisimple blocks, we infer  that $F$ preserves non-semisimple blocks.

Suppose that $F$ induces a stable equivalence between non-semisimple blocks $A_i$ and $B_j$ with $i,j\in [a_4]$. Then $U_i$ and $V_j$ are stably equivalent by Lemma \ref{sta-down}. In particular, $U_i$ and $V_j$ have the same number of non-projective indecomposable modules, that is, $n_i-1=m_j-1$. Hence $n_i=m_j$. $\square$

\medskip
Recall that the Auslander--Reiten conjecture on stable equivalences says that stably equivalent Artin algebras have the same number of non-projective simple modules (\cite[Conjecture (5), p.409]{ARS} and \cite{rouq}).

As a by-product of our methods developed in this subsection, we show the following result which slightly generalizes \cite[Theorem 1.4(1)]{xz3} that
states that the Auslander-Reiten conjecture holds true for stable equivalences between centralizer matrix algebras. Here, we do not assume that both algebras considered are centralizer matrix algebras.

\begin{Prop}\label{A-R}
Any centralizer matrix algebra over a field $R$ and its stably equivalent algebras have the same number of non-projective simple modules.
\end{Prop}

{\it Proof.} Let $c\in M_n(R)$ and $A:=S_n(c,R)$. We may assume $A=\prod^{a_1}_{i=1}A_i\times \prod^{a_2}_{i=a_1+1}A_i\times \prod^{l_c}_{i=a_2+1}A_i$ such that  $\{A_1,A_2,\cdots, A_{a_2}\}$ is exactly the set of all non-semisimple blocks of $A$, $\{A_1,A_2,\cdots, A_{a_1}\}$ is exactly the set of those non-semisimple blocks without nodes, and $\{A_{a_2+1}, \cdots, A_{l_c}\}$ is the set of simple blocks of $A$. Let $\widetilde{A}:=\prod^{a_1}_{i=1}A_i\times \prod^{a_2}_{i=a_1+1}(A_i)'$. Clearly, $\widetilde{A}$ has neither nodes nor semisimple direct summands. By Lemma \ref{node}(2), we have $n_i=LL(Z(A_i))\geq 3$ for $i\in [a_1]$, and $n_i=LL(Z(A_i))=2$ for $a_1<i\leq a_2$.

Now, suppose that an $R$-algebra $B$ is stably equivalent to $A$. Then the non-semisimple parts of $A$ and $B$ are stably equivalent. Let $\widetilde{B}$ be the product of non-semisimple blocks of the node-eliminated algebra $B'$ of $B$. Then $\widetilde{B}$ has neither nodes nor semisimple direct summands. Moreover, $\widetilde{A}$ and $\widetilde{B}$ are stably equivalent.
By the definition of $\widetilde{A}$ and $\widetilde{B}$, we see that $A$ (respectively, $B$) and $\widetilde{A}$ (respectively, $B$) have the same number of non-projective simple modules.

We shall show that $\widetilde{A}$ and $\widetilde{B}$ have the same number of non-projective simple modules. This implies that $A$ and $B$ have the same number of non-projective simple modules.

Indeed, according to Lemma \ref{exa}(2), the Frobenius parts of $\widetilde{A}$ and $\widetilde{B}$ are also stably equivalent. Further, the Frobenius part of $\widetilde{A}$ is Morita equivalent to $\prod^{a_1}_{i=1}U_i$ by Lemma \ref{node}(3) and Remark \ref{rmk2.11}, and the latter is a Nakayama algebra. Thanks to \cite[Theorem 1.3]{IR2}, an Nakayama Artin algebra and its stably equivalent Artin algebras have the same number of non-projective simples. Thus the Frobenius parts of $\widetilde{A}$ and $\widetilde{B}$ have the same number of non-projective simple modules. Now, according to \cite[Lemma 5.1(2)]{xz3}, we conclude that $\widetilde{A}$ and $\widetilde{B}$ have the same number of non-projective simple modules. $\square$

\medskip
Suppose that a non-semisimple block $A_i$ of $S_n(c,R)$ is stably equivalent to a non-semisimple block $B_j$ of $S_m(d,R)$. Then $U_i=R[x]/(f_i(x)^{n_i})$ and $V_j=R[x]/(g_j(x)^{m_j})$ are stably equivalent by Lemma \ref{sta-down}. Suppose further that $f_i(x)$ and $g_j(x)$ are separable. Then $ U_i\simeq V_j$ as $R$-algebras by Lemma \ref{poi}. Thus $N_j$ can be regarded as a $U_i$-module via this isomorphism, and is actually a generator for $U_i\modcat$. Subsequently, we need to study stable equivalences between the endomorphism algebras of generators for $A\modcat$, where $A$ is of the form $R[x]/(f(x)^n)$ for an irreducible polynomial $f(x)\in R[x]$ and $n\ge 2$.

\subsection{Stable equivalences between endomorphism algebras of generators}
In this subsection, we characterize stable equivalences between the endomorphism algebras of generators over quotients of the polynomial algebra in one variable.

Throughout this subsection, we fix an irreducible polynomial $f(x)\in R[x]$ of positive degree $u$ and set $A:=R[x]/(f(x)^n)$ for a positive integer $n\geq 2$. Let $M(i):= R[x]/(f(x)^i)\in A\modcat$ for $1\le i\le n$. We set $M(0):=0$. For $0\le j\le i \le n$, we denote by
$$f_{ji}: M(j)\ra M(i), \; x^k+(f(x)^j)\mapsto f(x)^{i-j}\cdot x^k+(f(x)^i) \mbox{ for } 0\le k\le uj-1 \mbox{ and } $$ $$g_{ij}: M(i)\ra M(j), \; x^k+(f(x)^i)\mapsto x^k+(f(x)^j) \mbox{ for } 0\le k\le ui-1$$ the canonical injective and surjective homomorphisms, respectively.
Clearly, $f_{ji}g_{ik}=0$ if $j+k \le i.$ For $j\le i$, there is an exact sequence
$0\to M(j)\to M(i)\to M(i-j)\to 0$ in $A\modcat$.

Any stable equivalence $G$ between $A$ and itself induces a permutation $\overline{G}$ on $\Gamma_{n-1}:=\{M(i)\mid i\in[n-1]\}$. For the stable equivalence $\overline{\Omega}_A: \stmc{A}\to A\stmc$ induced by the syzygy operator $\Omega_A$, we have $\overline{\Omega}_A^2=id.$

\smallskip
We start with consideration of the endomorphism algebras of generators for $A\modcat$.

\begin{Lem}\label{opp}
$(1)$ For $M\in A\modcat$, we have $\End_A(A\oplus M)\simeq \End_A(A\oplus M)\opp$ as $R$-algebras.

$(2)$ For $1\le j< i\le n$, let $\Lambda:=\End_A(A\oplus M(i))$. Then the exact sequence
$\eta_j: 0\ra M(j)\stackrel{f_{ji}}\ra M(i)\stackrel{g_{i i-j}}\ra M(i-j)\ra 0$ of $A$-modules induces an exact sequences of $\Lambda$-modules:
$$ 0\lra \Hom_{\Lambda}(A\oplus M(i), M(j))\stackrel{(f_{ji})_*}\lra \Hom_{\Lambda}(A\oplus M(i),M(i))\stackrel{(g_{i i-j})_*}\lra \Hom_{\Lambda}(A\oplus M(i),M(i-j))\lra 0.$$
\end{Lem}

{\it Proof.} (1) This  is an immediate consequence of Lemma \ref{iso-pr}.

(2) Note that the sequence $\Hom_A(A,\eta_j)$ is exact. We need to show that $\Hom_A(M(i),\eta_j)$ is also an exact sequence. Since $j<i$, the modules $M(i), M(j)$ and $M(i-j)$ can be viewed as $B:=R[x]/(f(x)^i)$-modules and we have $\Hom_A(M(i),M(k))=\Hom_B(M(i),M(k))$ for $k=i,j$ or $i-j$. As the $B$-module $M(i)$ is projective, $\Hom_B(M(i),\eta_j)$ is an exact sequence. Thus it follows from $\Hom_A(M(i), M(k))=\Hom_B(M(i), M(k))$ for $k=i,j$ or $i-j$ that $\Hom_A(M(i),\eta_j)$ is an exact sequence. $\square$

\smallskip
The following lemma will be used in characterizing stable equivalences between the endomorphism algebras of generators for $A\modcat$.

\begin{Lem}\label{endo}
Let $1\le j\leq i\le n-1$ and $\Lambda:=\End_A(A\oplus M(i))$. Then  $\End_{\Lambda}(\Hom_A(A\oplus M(i), M(j)))\simeq \End_A(M(j))=R[x]/(f(x)^j)$ as $R$-algebras.
\end{Lem}

{\it Proof.} Let $M:=A\oplus M(i)$ and $F:=\Hom_A(M,-)$ be the Hom-functor: $ A\modcat \ra \Lambda\modcat$. Then $$\End_{\Lambda}\big(F(M(j))\big)=\Hom_{\Lambda}\big(\Hom_A(M,M(j)),\Hom_A(M, M(j))\big)\simeq \Hom_A\big(_AM\otimes_{\Lambda}\Hom_A(M, M(j)), M(j)\big).$$
By Lemma \ref{opp}(2), the exact sequence $0\ra M(i-j)\ra M(i)\stackrel{g_{i j}}{\ra} M(j)\ra 0$ induces an exact sequence
$$ 0\lra F(M(i-j))\lra F(M(i))\stackrel{(g_{ij})_*}{\lra} F(M(j))\lra 0.$$
This shows that $g_{ij}$ is a right $\add(_AM)$-approximation of $M(j)$. Clearly, there is a surjective homomorphism $f: M\ra M(i-j)$. Thus $M(j)\in \app(M)$. By \cite[Lemm 2.2(1)]{x}, we have
$\End_{\Lambda}(F(M(j)))\simeq \End_A(M(j))$ as $R$-algebras.
$\square$

\medskip
The following lemma is useful in later proofs.

\begin{Lem}\label{indecom}
Let $j\in [n]$, $T:=\Hom_A(A\oplus M(j), A)$, $\Lambda:=\End_A(A\oplus M(j))$ and $T^*:=\Hom_{\Lambda}(T,\Lambda)$. Then the $\Lambda^{\opp}$-module $D\big(\Hom_A(A\oplus M(j), M(i))\big)$ is an indecomposable object in $\pres(T^*)$ for $i\in [j]$, where $D:=\Hom_R(-,R)$ is the $R$-duality. Moreover, $D\big(\Hom_A(A\oplus M(j), M(i))\big)\simeq{} _{\Lambda\opp}T^*\otimes_{\End_{\Lambda\opp}(T^*)} M(i)$ as $\Lambda\opp$-modules.
\end{Lem}

{\it Proof.} Let $H$ be the Hom-functor $\Hom_A(A\oplus M(j),-):A\modcat\to \Lambda\modcat$. Note that $T\simeq \nu_\Lambda T$ by our proof in Lemma \ref{node}(3). Therefore $T^*\simeq D(T)$ holds as $\Lambda^{\opp}$-modules. Let $\ires(_{\Lambda}T)$ be the full subcategory of $\Lambda\modcat$ consisting of those modules $L$ such that there is a minimal injective presentation $0\ra {}_{\Lambda}L\ra I_1\ra I_2$ with $I_1,I_2\in \add(_{\Lambda}T)$. Then $D:\ires(_{\Lambda}T)\ra \pres(T^*)$ is a duality of categories. By Lemma \ref{endo}, the modules $H(M(i))$, $i\in [j]$, are pairwise non-isomorphic indecomposable $\Lambda$-modules. Thus, to complete the proof, it suffices to show $H(M(i))\in \ires(_{\Lambda}T)$ for all $i\in [j]$.

Since $H$ is left exact, we may assume that $H(M(1))\subseteq H(M(2))\subseteq \cdots \subseteq P:=H(M(j))\subseteq T$ is a chain of submodules of $_{\Lambda}T$. Thus $\soc(_{\Lambda}H(M(k))\subseteq \soc(_{\Lambda}T)$ for $k\in[j]$. Due to $(T^*)_{\Lambda}\simeq D(T)_{\Lambda}$, we have $\soc(_{\Lambda}T)\simeq \top(_{\Lambda}T)$. By Lemma \ref{opp}(2), $P/H(M(i))\simeq H(M(j-i))$ for $i<j$. This implies  $\soc(P/H(M(i)))\simeq \soc(H(M(j-i)))\in \add\big(\soc(_{\Lambda}T)\big)$. To prove $H(M(i))\in \ires(_{\Lambda}T)$ for all $i\in [j]$, it suffices to show $\soc(T/H(M(i)))\in \add\big(\soc(_{\Lambda}T)\big)$ for all $i\in [j]$.

Suppose \emph{contrarily} that  $\soc(T/H(M(i)))$ has a simple submodule, say $V/H(M(i))$, such that $V/H(M(i))$ $\not\simeq \soc(_{\Lambda}T)$, where $V$ is a submodule of $T$ containing $H(M(i))$. On the one hand, as $\top(_{\Lambda}P)$ and $\soc(_{\Lambda}T)$ are all the simple $\Lambda$-modules (up to isomorphism), we have $ V/H(M(i))\simeq \top(_{\Lambda}P)$. From $\soc(P/H(M(i)))\in \add\big(\soc(_{\Lambda}T)\big)$,we have $\big(V/H(M(i))\big)\cap \big(P/H(M(i))\big)=0$. Thus $\top(_{\Lambda}P)\simeq V/H(M(i))$ is a composition factor of $(T/H(M(i)))/(P/H(M(i)))\simeq T/P$. On the other hand, the multiplicity $[L:S]$ of a simple module $S$ as composition factors of a module $L$ over an algebra $\Delta$ is given by $$[L:S]=\ell\big(_{\End_\Delta(P(S))}\Hom_\Delta(P(S),L)\big),$$ where $P(S)$ is a projective cover of $S$. In our case,  $\End_\Lambda(P)\simeq \End_A(M(j))\simeq R[x]/(f(x)^j)$, which is a local algebra with the unique simple module $R[x]/(f(x))$. Thus, for an $\End_\Lambda(P)$-module $X$, we have $\dim_R(X)=[X: \top(P)]\dim_R\big(R[x]/(f(x))\big)$, that is, the composition length of any $\End_\Lambda(P)$-module is determined by its $R$-dimension. Since $A$ is a symmetric algebra, we have
$\dim_R(\Hom_\Lambda(P,T))=\dim_R(\Hom_A(M(j),A))=\dim_R(M(j)).$ Further,  $\dim_R(\Hom_\Lambda(P,P))=\dim_R(\Hom_A(M(j),M(j)))=\dim_R(M(j)).$
Hence the multiplicities of $\top(_{\Lambda}P)$ as composition factors in both $P$ and $T$ are equal. Thus $T/P$ cannot have $\top(_{\Lambda}P)$ as a composition factor. This proves the first statement.

Now, we prove the second statement. Since $T^*\simeq D(T)$ as $\Lambda^{\opp}$-modules, it follows that $$\End_{\Lambda\opp}(T^*)\simeq \End_{\Lambda\opp}(D(T))\simeq \End_\Lambda(T)\opp\simeq A.$$ So we can identify $\End_{\Lambda\opp}(T^*)$ with $A$. As $T^*\otimes_A -: A\modcat \ra \pres(T^*)$ is an equivalence of categories, the functor $T^*\otimes_A -$ preserves indecomposable objects and endomorphism algebras of modules. Note that two indecomposable $A$-modules $X$ and $Y$ are isomorphic if and only if $\End_A(X)\simeq \End_A(Y)$ as algebras. Consequently, two indecomposable objects in $\pres(T^*)$ are isomorphic if and only if their endomorphism algebras are isomorphic. Thus it is sufficient to show that $DH(M(i))$ and $ _{\Lambda\opp}T^*\otimes_A M(i)$ have isomorphic endomorphism algebras. However, this is clear from the following isomorphisms:
$$\End_{\Lambda\opp}(DH(M(i)))\simeq \End_\Lambda(H(M(i)))\opp\simeq \End_A(M(i))\opp\simeq R[x]/(f(x)^i)$$ and
$ \End_{\Lambda\opp}(_{\Lambda\opp}T^*\otimes_A M(i))\simeq \End_A(M_i)\simeq R[x]/(f(x)^i).$ $\square$

\medskip
Now, let $\Lambda$ be an arbitrary algebra. For $P\in \prj{\Lambda}$, we denote by $\underline{{\rm pres}}(P)$ the full subcategory of $\Lambda\stmc$ with objects in $\pres(P)$ (see Section \ref{sect2.1} for definition).

Suppose that there is a stable equivalence $F:\Lambda\stmc{}\ra \Gamma\stmc$ between algebras $\Lambda$ and $\Gamma$. Let $T\in \prj{\Lambda}$, $W\in \prj{\Gamma}$, $E:=\End_\Lambda(T)$ and $L:=\End_\Gamma(W)$. If $F$ restricts to an equivalence $F_T: \underline{{\rm pres}}(T)\ra \underline{{\rm pres}}(W)$, then there exists the following commutative diagram of functors (up to natural isomorphism):

$$\xymatrix@C=1.5cm@R=0.8cm{
&\Lambda\stmc{}\ar[r]^-{F}  & \Gamma\stmc{}\\
(\epsilon_{(F_{\Lambda, \Gamma}, T, W)}): &\underline{{\rm pres}}(T)\ar@{^{(}->}[u]\ar[r]^-{F_T} & \underline{{\rm pres}}(W) \ar@{^{(}->}[u]\\
&E\stmc{}\ar[r]^-{G}\ar[u]^-{T\otimes_E-}  & \Gamma\stmc{}\ar[u]_-{W\otimes_E-}\\
}
$$
with $G:=\Hom_\Gamma(W,-)\circ F_T\circ (T\otimes_E -)$ being an equivalence, where $T\otimes_E -:E\stmc{}\ra \underline{{\rm pres}}(T) $ and $W\otimes_L -:L\stmc{}\ra \underline{{\rm pres}}(W)$ are the equivalences induced by the equivalences $T\otimes_E -: E\modcat \ra \pres(T)$ and $W\otimes_L - : L\modcat\ra \pres(W)$, respectively. Observe that $T\otimes_E -:E\stmc{}\ra \underline{{\rm pres}}(T) $ is well defined because a homomorphism $f: X\to Y$ of modules factorizes through a projective module if and only if $f$ factorizes through a projective cover $P(Y)$ of $Y$, while $P(Y)$ lies in $\add(T)$ if $Y\in {\rm pres}(T)$.

\begin{Def} Let $T\in \Lambda\prj$. A stable equivalence $F:\Lambda\stmc{}\ra \Gamma\stmc$ of Artin algebras $\Lambda$ and $\Gamma$ satisfies $(\delta_T)$ if $\End_\Lambda(X)\simeq \End_\Gamma(F(X))$ for all non-projective indecomposable $\Lambda$-modules $X\in \pres(T)$.
\end{Def}

\begin{Lem}\label{end}
Given $(\epsilon_{(F_{\Lambda, \Gamma}, T, W)})$ and $X\in E\prj$, the functor $G$ in the diagram satisfies $(\delta_X)$ if and only if $F$ satisfies $(\delta_{T\otimes_E X})$.
\end{Lem}

{\it Proof.} Clearly, $T\otimes_E X\in \add(T)$ and therefore is a projective $\Lambda$-module. Recall that $G$ induces a one-to-one correspondence $G: E\modcat_{\mathscr{P}}\ra L\modcat_{\mathscr{P}}$, which clearly preserves non-projective indecomposable modules. Let $Z$ be an arbitrary non-projective indecomposable $E$-modules in $\pres(X)$. Then $G(Z)$ is a non-projective indecomposable $L$-module. Since $T\otimes_E -: E\modcat \ra \pres(T)$ is an equivalence of categories, it follows that $T\otimes_E Z$ is a non-projective indecomposable module in $\pres(T\otimes_E X)$ and $\End_E(Z)\simeq \End_\Lambda(T\otimes_E Z)$ as algebras. Similarly,  $W\otimes_L G(Z)$ is a non-projective indecomposable modules in $\pres(W)$ and $\End_L(G(Z))\simeq \End_\Gamma(W\otimes_L G(Z))$ as algebras. It follows from the commutative diagram $(\epsilon_{(F_{\Lambda, \Gamma}, T, W)})$ that $F(T\otimes_E Z)=F_T(T\otimes_E Z)\simeq W\otimes_L G(Z)$ in $\Gamma\stmc$. Since $F(T\otimes_E Z)$ and $ W\otimes_L G(Z)$ are non-projective indecomposable modules, we have $F(T\otimes_E Z)\simeq W\otimes_L G(Z)$ in $\Gamma\modcat$. As $Z$ is arbitrary, we deduce that $G$ satisfies $(\delta_X)$ if and only if $\End_\Lambda(T\otimes_E Z)\simeq \End_\Gamma(W\otimes_L G(Z))\simeq \End_\Gamma(F(T\otimes_E Z))$ for all non-projective indecomposable module $Z\in \pres(X)$. Since $T\otimes_E -: E\modcat \ra \pres(T)$ is an equivalence, any indecomposable module $Y$ in $\pres(T\otimes_E X)$ is isomorphic to $T\otimes_E Z$ for some indecomposable module $Z$ in $\pres(X)$. Thus $\End_\Lambda(T\otimes_E Z)\simeq \End_\Gamma(F(T\otimes_E Z))$ for all non-projective indecomposable module $Z\in \pres(X)$ if and only if $F$ satisfies $(\delta_{T\otimes_E X})$. $\square$

\medskip
\textbf{In the rest of this subsection}, we assume that $M, N\in A\modcat_{\mathscr{P}}$, $\Lambda:=\End_A(A\oplus M)$, $\Gamma:=\End_A(A\oplus N)$, $_{\Lambda}T:=\Hom_A(A\oplus M, A)$ and $_{\Gamma}W:=\Hom_A(A\oplus N, A)$. Then

(i) $\Lambda\simeq \Lambda^{\opp}$ and  $\Gamma\simeq \Gamma^{\opp}$ by Lemma \ref{opp}(1).

(ii) $T\simeq D(T^*)=\nu_\Lambda(T)$ (respectively, $W\simeq D(W^*)=\nu_\Lambda(W)$) is the unique indecomposable projective-injective $\Lambda$-module (respectively, $\Gamma$-module) by the proof of Lemma \ref{node}(3).

(iii) $\End_\Lambda(T)\simeq A\simeq \End_{\Lambda\opp}(T^*)$ and $\End_\Gamma(W)\simeq A\simeq \End_{\Gamma\opp}(W^*)$.

The purpose of the rest of this subsection is to establish a relation between $M$ and $N$ under a stable equivalence between $\Lambda$ and $\Gamma$.

\begin{Lem}\label{Stc-pre} Suppose that $n\geq 3$, $M$ and $N$ are indecomposable. If $F: \Lambda\stmc{}\ra \Gamma\stmc$ is a stable equivalence, then there is a commutative diagram $(\epsilon_{(F_{\Lambda\opp,\Gamma\opp}, T^*,W^*)})$. Further, if $\overline{G}$ is the bijection induced by the stable equivalence $G:A\stmc{}\ra A\stmc$ in this diagram, then $\overline{G}=id$ if and only if $F$ satisfies $(\delta_T)$.
\end{Lem}

{\it Proof.}
Let $H_M$ denote the Hom-functor $\Hom_A(A\oplus M,-):A\modcat\ra \Lambda\modcat$, $T:=H_M(A)$ and $W:=H_N(A)$. Then $T$ (respectively, $W$) is the unique projective-injective indecomposable $\Lambda$-module (respectively, $\Gamma$-module). It follows from $n\geq 3$ and Lemma \ref{node} that both $\Lambda$ and $\Gamma$ have no nodes. Clearly, $\Lambda$ and $\Gamma$ have no nonzero semisimple summands. According to Lemma \ref{exa1}, the stable equivalence functor $F$ induces an equivalence $\underline{\pres}(T)\ra \underline{{\rm pres}}(W)$. Thus we have a commutative diagram $(\epsilon_{(F_{\Lambda,\Gamma}, T, W)})$.

\smallskip
Since $T\simeq \nu_\Lambda(T)$ as $\Lambda$-modules and $W\simeq \nu_\Gamma(W)$ as $\Gamma$-modules, it follows that $T^*:=\Hom_\Lambda(T,\Lambda)\simeq D(T)$ and $W^*:=\Hom_\Gamma(W,\Gamma)\simeq D(W)$ are the unique projective-injective indecomposable modules over $\Lambda\opp$ and $\Gamma\opp$, respectively. By (i), %Lemma \ref{opp},
we have $\Lambda\simeq \Lambda\opp$ and $\Gamma\simeq \Gamma\opp$. So, if $\Lambda\modcat$ (respectively, $\Gamma\modcat$) is identified with $\Lambda\opp\modcat$ (respectively, $\Gamma\opp\modcat$) via those isomorphisms, then $T$ is identified with $T^*$ and $W$ is identified with $W^*$, and therefore $\pres(T)=\pres(T^*)$ and $\pres(W)=\pres(W^*)$. Thus we obtain the desired commutative diagram $(\epsilon_{(F_{\Lambda\opp,\Gamma\opp}, T^*,W^*)})$. This provides a stable equivalence $G:\stmc{A}\to A\stmc$.

According to Lemma \ref{self}, the permutation $\overline{G}$, induced by $G$, on the set $\Gamma_{n-1}$ equals either $id$ or $\overline{\Omega}_A$.  It follows from $n\geq 3$ that $\End_A(M(1))\not\simeq\End_A(M(n-1))$. Thus $\overline{G}=id$ if and only if $\End_A(M(i))\simeq \End_A(\overline{G}(M(i)))$ as algebras for all $i\in [n-1]$. By Lemma \ref{end}, we deduce that $\End_A(M(i))$ $\simeq $ $\End_A(\overline{G}(M(i)))$ as algebras for all $i\in [n-1]$ if and only if $F$ satisfies $(\delta_{T^*})$. Since $\pres(T)=\pres(T^*)$ by identifying $\Lambda$ with $\Lambda\opp$ (respectively, $\Gamma$ with $\Gamma\opp$), we see that $F$ satisfies $(\delta_{T^*})$ if and only if $F$ satisfies $(\delta_T)$. $\square$

\begin{Lem}\label{Stc} Suppose that $M$ and $N$ are indecomposable and that $F: \Lambda\stmc{}\ra \Gamma\stmc$ is a stable equivalence. Then $_AM\simeq {}_AN$ if $F$ satisfies $(\delta_T)$, and $_AM\simeq \Omega_A(N)$, otherwise.
\end{Lem}

{\it Proof.} If $n=2$, then $A$ has $2$ indecomposable modules. Thus the statement follows. Now, assume $n\geq 3$. In this case, $\End_A(M(1))\not\simeq\End_A(M(n-1))$. Let $a=\ell(M)$ and $b=\ell(N)$. Then $M\simeq M(a)$ and $N\simeq M(b)$. Without loss of generality, assume $a\geq b$.

Suppose that $F:\Lambda\stmc{}\ra \Gamma\stmc$ is a stable equivalence. Let $H_M:=\Hom_A(A\oplus M,-):A\modcat\ra \Lambda\modcat$, $_{\Lambda}T:=H_M(A)\in \Lambda$-proj and $_{\Gamma}W:=H_N(A)\in \Gamma$-proj. Then $\End_{\Lambda^{\opp}}(T^*)\simeq A\simeq  \End_{\Gamma^{\opp}}(L^*)$, there are a diagram $(\epsilon_{(F_{\Lambda\opp,\Gamma\opp}, T^*,W^*)})$ and a stable equivalence $G:\stmc{A}\to A\stmc$ by Lemma \ref{Stc-pre}.

Set $\Delta:=\End_A(A\oplus \Omega_A(N))$. Then there is the stable equivalence $J: \Gamma\stmc{}\ra \Delta\stmc$ constructed in \cite[Section 3]{LX3}. According to Lemma \ref{Stc-pre}, $J$ induces a stable equivalence $K: A\stmc{}\ra A\stmc$. By the construction of $J$, the permutation $\overline{K}$, induced by $K$, on $\Gamma_{n-1}$ equals $\overline{\Omega}_A$. According to Lemma \ref{self}, the permutation $\overline{G}$, induced by $G$, on the set $\Gamma_{n-1}$ equals either $id$ or $\overline{\Omega}_A$. Hence either $\overline{G}=id$ or $\overline{KG}=\overline{K}\;\overline{G}=\overline{\Omega}_A^2=id$. By Lemma \ref{Stc-pre}, $\overline{G}=id$ if and only if $F$ satisfies $(\delta_T)$. Thus, to complete the proof, it suffices to show that $M\simeq N$ when $\overline{G}=id$.

Suppose $\overline{G}=id$. To prove $M\simeq N$, we only need to show $a=b$.

By Lemma \ref{indecom}, $T^*\otimes_A M(k) \simeq D(H_M(M(k))$ in $\Lambda\opp\modcat$ and $W^*\otimes_A M(k) \simeq D(H_N(M(k))$ in $\Gamma\opp\modcat$ for $k\in [b]$.
From the commutative diagram $(\epsilon_{(F_{\Lambda\opp,\Gamma\opp}, T^*,W^*)})$ and the assumption $\overline{G}=id$ it follows that
$$F(D(H_M(M(k))))\simeq F(T^*\otimes_A M(k))\simeq W^*\otimes_A G(M(k))\simeq W^*\otimes_A M(k)\simeq D(H_N(M(k)))$$
in $\Gamma\opp\stmc$ for $k\in [b]$. Note that two non-projective indecomposable $\Gamma\opp$-modules are isomorphic in $\Gamma\opp\stmc$ if and only if they are isomorphic as $\Gamma\opp$-modules. Thus $F(D(H_M(M(k))))\simeq D(H_N(M(k)))$ in $\Gamma\opp$-mod for $k\in [b]$.

Contrarily, suppose $a>b$.  We consider the following exact sequence of $A$-modules
$$\diamondsuit: \quad 0\lra M(a-b)\stackrel{l}\lra M(a)\stackrel{h}\lra M(b)\lra 0$$ with $l:=f_{a-b,b}$ and $h:=g_{a,b}$ the canonical homomorphisms, and apply the functor $DH_M$ to $(\diamondsuit)$. This gives rise to the following exact sequence of $\Lambda\opp$-modules
$$DH_M(\diamondsuit): \; 0\lra DH_M(M(b))\stackrel{DH_M(h)}\lra D(P)\stackrel{DH_M(l)}\lra DH_M(M(a-b))\lra 0.$$ Note that the exact sequence $DH_M(\diamondsuit)$ has no nonzero split direct summands. Applying Lemma \ref{exa}(3), we get the exact sequence of $\Gamma\opp$-modules:
$$DH_M(\diamondsuit)': 0\lra F(DH_M(M(b)))\stackrel{DH_M(h)'}\lra (W^*)^{\oplus m}\oplus F(D(P))\stackrel{DH_M(l)'}\lra F(DH_M(M(a-b)))\lra 0,$$where $m$ is a non-negative integer and $DH_M(\diamondsuit)'$ does not have nonzero split direct summands. But $F(DH_M(M(b)))\simeq DH_N(M(b)) = D(Q)$ is an injective $\Gamma\opp$-module, which yields that $DH_M(\diamondsuit)'$ does have a nonzero split direct summand. This is a contradiction and shows $a=b$. $\square$

\begin{Prop}\label{Stc2}
$\Lambda$ and $\Gamma$ are stably equivalent if and only if $\mathcal{B}(M)\simeq \mathcal{B}(N)$ or $\mathcal{B}(M)\simeq \Omega_A(\mathcal{B}(N)).$
\end{Prop}
{\it Proof.} If $n=2$, then the statement follows immediately. Now, assume $n\geq 3$. Without loss of generality, we assume that $M$ and $N$ are basic modules.

If $M\simeq N$, then $\Lambda\simeq \Gamma$, and therefore $\Lambda$ and $\Gamma$ are stably equivalent. If $M\simeq \Omega_A(N)$, then $\Lambda$ and $\Gamma$ are stably equivalent by Lemma \ref{almm}.

Now, suppose that $F:\Lambda\stmc{}\ra \Gamma\stmc$ is a stable equivalence. We shall prove either $M\simeq N$ or $M\simeq \Omega_A(N)$.

Indeed, it follows from $n\geq 3$ and Lemma \ref{node} that $\Lambda$ and $\Gamma$ have no nodes. According to Lemma \ref{exa}(1), $F$ induces a one-to-one correspondence $F':\mathscr{P}(\Lambda)_{\mathscr{I}}\ra \mathscr{P}(\Gamma)_{\mathscr{I}}$. Now we write $_AM=\oplus^r_{i=1}X_i$ and $_AN=\oplus^s_{j=1}Y_j$, where $X_i$ and $Y_j$ are indecomposable $A$-modules such that $X_i\not\simeq X_j$ for $i\neq j$ and $Y_p\not\simeq Y_q$ for $p\ne q$. Let $H_M:=\Hom_A(A\oplus M,-): A\modcat\ra \Lambda\modcat $, $H_N:=\Hom_A(A\oplus N, -): A\modcat\ra \Gamma\modcat$. Then $T:=H_M(A)$ and $W:=H_N(A)$ are the unique indecomposable projective-injective modules over $\Lambda$ and $\Gamma$, respectively, and the sets $\{H_M(X_i)\mid i\in [r]\}$ and $\{H_N(Y_j)\mid j\in [s]\}$ are complete representatives of isomorphism classes of non-injective, indecomposable projective modules over $\Lambda$ and $\Gamma$, respectively. The bijection of $F'$ implies $r=s$. Thus we may assume $F'((H_M(X_i))=H_N(Y_i)$ for $i\in [r]$.

For $i\in [r]$, set $T_i:=H_M(A\oplus X_i)$, $W_i:=H_N(A\oplus Y_i)$, $\Lambda_i:=\End_\Lambda(T_i)$ and $\Gamma_i:=\End_\Gamma(W_i)$. Then $T_i=T\oplus H_M(X_i)\in \Lambda\prj$ and $W_i=W\oplus H_N(Y_i)\in \Gamma\prj$.
If $Z\in \pres(T_i)$ is a non-projective indecomposable $\Lambda$-module, then $F(Z)\in \pres(W_i)$ by Lemma \ref{exa1}. Thus $F$ restricts to an $R$-linear equivalence $F_{T_i}: \underline{{\rm pres}}(T_i)\ra \underline{{\rm pres}}(W_i)$. As $\Lambda_i\modcat$ and ${\rm pres}(T_i)$ are equivalent, we have not only a stable equivalence $G_i:\stmc{\Lambda_i}\to \Gamma_i\stmc $, but also a commutative diagram $(\epsilon_{(F_{\Lambda,\Gamma}, T_i,W_i)})$.

Let $X'_i:=\Hom_\Lambda(T_i,T)\in \Lambda_i\prj$.
Due to $_{\Lambda}T\in \add(T_i)$, we have $T_i\otimes_{\Lambda_i}X'_i=T_i\otimes_{\Lambda_i} \Hom_{\Lambda}(T_i,T)\simeq T$ as $\Lambda$-modules, and therefore $\pres(T_i\otimes_{\Lambda_i}X'_i)=\pres(T)$. Applying Lemma \ref{end} to $(\epsilon_{F_{\Lambda,\Gamma},T_i, W_i})$, we deduce that $G_i$ satisfies $(\delta_{X'_i})$ if and only if $F$ satisfies $(\delta_{T})$. Hence, by Lemma \ref{Stc}, for  $i\in [r]$, we have $_AX_i\simeq {}_AY_i$ if $F$ satisfies $(\delta_{T})$, and $X_i\simeq \Omega_A(Y_i)$ otherwise. Since $i$ is arbitrary in $[r]$ and $T$ does not depend upon $i$, we get $M\simeq N$ or $M\simeq \Omega_A(N)$.
$\square$

\subsection{Stable equivalences of centralizer matrix algebras: Proofs of main results}

This subsection is devoted to the proofs of Theorem \ref{main1}.

\medskip
{\bf Proof of Theorem \ref{main1}.}
(1) Suppose that $S_n(c,R)$ and $S_m(d,R)$ are stably equivalent. By Lemma \ref{iso-pr}, $S_n(c,R)$ is semisimple if and only if $\mathcal{M}_c$ consists of only irreducible polynomials if and only if $\mathcal{R}_c=\emptyset$. Since an algebra stably equivalent to a semisimple algebra itself is semisimple, we may assume that both $S_n(c,R)$ and $S_m(d,R)$ are non-semisimple, that is, $\mathcal{R}_c\neq \emptyset$ and $\mathcal{R}_d\neq\emptyset$.

Let $A_1, \cdots, A_s$ and $B_1,\cdots, B_t$ be the non-semisimple blocks of $S_n(c,R)$ and $S_m(d,R)$, respectively, and  $A:=A_1\times\cdots \times A_s$ and $B:=B_1\times\cdots\times B_t$.

By Proposition \ref{sta}, we have $s=t$. For $i\in [s]$, the block $A_i$ is stably equivalent to $B_j$ for some $j\in [t]$. In this case, $n_i=m_j$. Thus we may assume that $A_i$ is stably equivalent to $B_i$ for $1\leq i\leq s$. Recall that $A_i:=\End_{U_i}(M_i)$ and $B_i:=\End_{V_i}(N_i)$ with $U_i:=R[x]/(f_i(x)^{n_i})$ and $V_i:=R[x]/(g_i(x)^{m_i})$ by our notation. Note that the non-semisimple blocks of $S_n(c,R)$ are in one-to-one correspondence to the elementary divisors in $\mathcal{R}_c$. Thus we may define a bijective map $\pi: \mathcal{R}_c\ra \mathcal{R}_d, f_i(x)^{n_i}\mapsto g_i(x)^{m_i}$ for $i\in [s]$. Note that $A_i$ is Morita equivalent to $\End_{U_i}(\mathcal{B}(M_i))$, and $B_i$ is Morita equivalent to $\End_{V_i}(\mathcal{B}(N_i))$. This implies that $\End_{U_i}(\mathcal{B}(M_i))$ and $\End_{V_i}(\mathcal{B}(N_i))$ are stably equivalent. Thus $U_i$ and $V_i$ are stably equivalent by Lemma \ref{sta-down}. By Lemma \ref{poi}(3), we get $U_i\simeq V_i$ as algebras, that is, $R[x]/(f_i(x)^{n_i})\simeq R[x]/(g_i(x)^{m_i})= R[x]/\big((f_i(x)^{n_i})\pi\big)$. Via this isomorphism, we may view $N_i$ as a module over $U_i$ for $i\in [s]$. Then $\mathcal{B}(N_i)\simeq \bigoplus_{j\in {P_d(g_i(x)^{m_i})}} R[x]/(f_i(x)^j)$ as $U_i$-modules. Recall that $\mathcal{B}(M_i)\simeq \bigoplus_{r\in {P_c(f_i(x)^{n_i})}} R[x]/(f_i(x)^r)$ as $U_i$-modules. Since $M_i$ and $N_i$ are faithful $U_i$-modules, respectively, it follows that $\mathcal{B}(M_i)\simeq U_i\oplus \mathcal{B}(M_i)_{\mathscr{P}}$ and $\mathcal{B}(N_i)\simeq U_i\oplus \mathcal{B}(N_i)_{\mathscr{P}}$. By Proposition \ref{Stc2}, either $\mathcal{B}(M_i)_{\mathscr{P}}\simeq \mathcal{B}(N_i)_{\mathscr{P}}$ or $\mathcal{B}(M_i)_{\mathscr{P}}\simeq\Omega_{U_i}(\mathcal{B}(N_i)_{\mathscr{P}})$. Clearly, we have ${P_c(f_i(x)^{n_i})}={P_d(g_i(x)^{m_i})}$ if $\mathcal{B}(M_i)_{\mathscr{P}}\simeq \mathcal{B}(N_i)_{\mathscr{P}}$, and $P_c(f_i(x)^{n_i})=\mathcal{J}_{P_d(g_i(x)^{m_i})}$ if $\mathcal{B}(M_i)_{\mathscr{P}}\simeq \Omega_{U_i}(\mathcal{B}(N_i)_{\mathscr{P}})$. Thus $c$ and $d$ are $S$-equivalent.

Conversely, suppose $c\stackrel{S}\sim d$. Then, by definition, there is a bijection $\pi$ between $\mathcal{R}_c$ and $\mathcal{R}_d$ such that, for $f_i(x)^{n_i}\in \mathcal{R}_c$, the isomorphism $R[x]/(f_i(x)^{n_i})\simeq R[x]/((f_i(x)^{n_i})\pi)$ holds as algebras and either ${P_c(f_i(x)^{n_i})}={P_d((f_i(x)^{n_i})\pi)}$ or $P_c(f_i(x)^{n_i})=\mathcal{J}_{P_d((f_i(x)^{n_i})\pi)}$. Thus $s=t$. 
We may assume $(f_i(x)^{n_i})\pi=g_i(x)^{m_i}$ for $i\in [s]$. In this case we have $n_i=m_i$. Then $U_i=R[x]/(f_i(x)^{n_i})\simeq R[x]/(g_i(x)^{m_i})=V_i$ for $i\in [s]$. This isomorphism enables us to view $N_i$ as a $U_i$-module. Now the condition ${P_c(f_i(x)^{n_i})}={P_d(g_i(x)^{m_i})}$ and $P_c(f_i(x)^{n_i})=\mathcal{J}_{P_d(g_i(x)^{m_i})}$ are equivalent to the isomorphisms $\mathcal{B}(M_i)_{\mathscr{P}}\simeq \mathcal{B}(N_i)_{\mathscr{P}}$ and $\mathcal{B}(M_i)_{\mathscr{P}}\simeq \Omega_{U_i}(\mathcal{B}(N_i)_{\mathscr{P}})$, respectively. Hence $A_i=\End_{U_i}(M_i)$ and $B_i=\End_{V_i}(N_i)$ are stably equivalent in both cases for $i\in [s]$ (by applying Lemma \ref{almm} to the case $\mathcal{B}(M_i)_{\mathscr{P}}\simeq \Omega_{U_i}(\mathcal{B}(N_i)_{\mathscr{P}})$). This yields that $A = \prod_i^sA_i$ and $B= \prod_i^sB_i$ are stably equivalent. Hence $S_n(c,R)$ and $S_m(d,R)$ are stably equivalent.

(2) This is proved in Proposition \ref{A-R}.
$\square$

\subsection{Stable equivalences for permutation matrices: Proof of Corollary \ref{cor1.4}\label{stper}}
This subsection is devoted to proving Corollary \ref{cor1.4}. We first prepare several useful lemmas.

\medskip
Let $\sigma=\sigma_1\cdots\sigma_s\in \Sigma_n$ be a product of disjoint cycle-permutations $\sigma_i$, and let $\lambda=(\lambda_1,\cdots, \lambda_s)$ be its cycle type with $\lambda_i\ge 1$ for $i\in [s]$. For a prime number $p>0$,  a cycle $\sigma_i$ is said to be \emph{$p$-regular} if $p\nmid \lambda_i$, and \emph{$p$-singular} if $p\mid \lambda_i$. If $p=0$, all cycles are regarded as $p$-regular cycles. Let $r(\sigma)$ (respectively, $s(\sigma)$) be the product of the $p$-regular (respectively, $p$-singular) cycles of $\sigma$. We consider $r(\sigma)$ and $s(\sigma)$ as elements in $\Sigma_n$ by letting $r(\sigma)$ (respectively, $s(\sigma)$) fix the points involved in the $p$-singular cycles (respectively, $p$-regular cycles) of $\sigma$. The permutation $\sigma$ is said to be \emph{$p$-regular} (or \emph{$p$-singular}) if  $\sigma = r(\sigma)$ (or $\sigma=s(\sigma)$).
For $p=0$, we define $r(\sigma)=\sigma$ and $s(\sigma)=id$, the identity permutation.

We denote by $c_{\sigma}:=\sum_{i=1}^ne_{i,(i)\sigma}\in M_n(R)$ the permutation matrix of $\sigma$ over $R$, where $e_{ij}$ is the matrix with $1$ in $(i,j)$-entry and $0$ in other entries.

Suppose that $g(x)$ is an irreducible factor of the minimal polynomial $m_{c_\sigma}(x)$ of $c_{\sigma}$. We define
$$q_{g(x)}:= max\{\nu_{p}(\lambda_j) \mid j\in[k], \mbox{ such that } g(x)\mbox{ divides } x^{\lambda_j}-1 \}.$$
The elementary divisors of permutation matrices can be described in terms of their cycle types.

\begin{Lem}{\rm \cite[Lemma 2.17]{lx1}}\label{per}
Suppose that the characteristic of $R$ is $p\ge 0$ and $\sigma\in \Sigma_n$ is a permutation of cycle type $(\lambda_1,\cdots,\lambda_k)$. Then $\mathcal{E}_{c_{\sigma}} = \{g(x)^{p^{\nu_p(\lambda_i)}}\mid i\in[k], \; g(x)~\mbox{is an irreducible factor of}~ x^{\lambda_i}-1 \}.$
\end{Lem}
Let $c$ be a permutation matrix over a field of characteristic $p\geq 0$. Then, by Lemma \ref{per}, the matrix $c$ always has a maximal elementary divisor of the form $(x-1)^{p^a}$ for some non-negative integer $a\geq 0$. It is uniquely determined by $c$, and is called the \emph{exceptional} elementary divisor.

\smallskip
The following result is easy,  we omit its proof.

\begin{Lem}\label{p-p}
Let $p>0$ be a prime, and let $S, T \subseteq \mathbb{Z}_{>0}$ be two sets consisting of $p$-powers. Suppose $S=\mathcal{J}_T$. Then either $S=T$ is a singleton set, or $p=2$ and $S=T=\{2^u,2^{u+1}\}$ for some $u\in\mathbb{N}$.
\end{Lem}
As a corollary of Lemma \ref{per}, we have the following.
\begin{Lem}\label{map}
Let $R$ be a field of characteristic $p>0$, and  let $\sigma\in \Sigma_n, \tau\in \Sigma_m$ be two permutations such that $p$ divides the orders of both $\sigma$ and $\tau$. If $c_\sigma\stackrel{S}\sim c_\tau$, then there exists a bijection $\pi': \mathcal{R}_{c_\sigma}\ra \mathcal{R}_{c_\tau}$ such that

$(1)$ $\pi'$ defines an $S$-equivalence between $c_\sigma$ and $c_\tau$.

$(2)$ $\pi'$ maps the exceptional elementary divisor of $c_\sigma$ to the exceptional elementary divisor of $c_\tau$.

\end{Lem}
{\it Proof.} Let $(x-1)^{p^a}$ and $(x-1)^{p^b}$ be the exceptional elementary divisors of $c_\sigma$ and $c_\tau$, respectively, where $a$ and $b$ are nonnegative integers. Since $p$ divides the order of  $\sigma$ (respectively,  $\tau$), it follows that $p$ divides at least one component of the cycle type of $c$ (respectively, $d$). Thus, by Lemma \ref{per}, we have $a>0$ and $b>0$. This implies that $(x-1)^{p^a}\in \mathcal{R}_{c_\sigma}$ and $(x-1)^{p^b}\in \mathcal{R}_{c_\tau}$.

Suppose $c_\sigma\stackrel{S}\sim c_\tau$. Then, by definition, there is a bijection $\pi: \mathcal{R}_{c_\sigma}\ra \mathcal{R}_{c_\tau}$  such that, for $f(x)\in \mathcal{R}_c$, the two conditions hold:

$(1)$ $R[x]/(f(x))\simeq R[x]/((f(x))\pi)$ as algebras, and

$(2)$ either $P_c(f(x))= P_d((f(x))\pi)$ or $P_c(f(x))=\mathcal{J}_{P_d((f(x))\pi)}$.

Clearly, it follows from Lemma \ref{per} that both $P_c(f(x))$ and $P_d((f(x))\pi)$ consist of $p$-powers for all $f(x)\in \mathcal{R}_c$. Thus, by Lemma \ref{p-p}, we deduce that $P_c(f(x))= P_d((f(x))\pi)$ for all $f(x)\in \mathcal{R}_c$. If $((x-1)^{p^a})\pi=(x-1)^{p^b}$, then $\pi':=\pi$ is a desired map.

Now, assume $((x-1)^{p^a})\pi\neq (x-1)^{p^b}$. For brevity, let $j(x):=((x-1)^{p^a})\pi$ and $k(x):=((x-1)^{p^b})\pi^{-1}$. Then it follows from  $R[x]/(h(x))\simeq R[x]/((h(x))\pi)$ for $h(x)\in \mathcal{R}_{c_\sigma}$ that $j(x)=(x+u)^{p^a}$ and $k(x)=(x+v)^{p^b}$ for $u,v\in R$ with $u\ne -1\ne v$. By definition of $\pi$, $P_{c_{\sigma}}((x-1)^{p^a})=P_{c_{\tau}}((x+u)^{p^a})$ and $P_{c_{\sigma}}((x+v)^{p^b})=P_{c_{\tau}}((x-1)^{p^b})$. Further, by Lemma \ref{per}, we have $P_{c_{\sigma}}((x+v)^{p^b})\subseteq P_{c_{\sigma}}((x-1)^{p^a})$ and $P_{c_{\tau}}((x+u)^{p^a})\subseteq P_{c_{\tau}}((x-1)^{p^b})$. Thus the two inclusions and the two equalities show the following:
$$P_{c_{\sigma}}((x+v)^{p^b})=P_{c_{\tau}}((x-1)^{p^b})= P_{c_{\sigma}}((x-1)^{p^a})=P_{c_{\tau}}((x+u)^{p^a}).$$ This implies that both $p^a$ and $p^b$ are the maximal numbers in both $P_{c_{\sigma}}((x-1)^{p^a})$ and $P_{c_{\tau}}((x-1)^{p^b})$. Thus $p^a=p^b$ and $a=b$. Note that $\mathcal{R}_{c_{\sigma}}=\mathcal{R}_{c_{\sigma}}\setminus\{(x-1)^{p^a}, (x-v)^{p^a}\}\cup \{(x-1)^{p^a}, (x-v)^{p^a}\}$. We define a map \vspace{-0.3cm}
$$ \pi': \mathcal{R}_{c_{\sigma}}\ra \mathcal{R}_{c_{\tau}}, \quad h(x)\mapsto \begin{cases}(h(x))\pi & \mbox{ if } h(x)\not\in \{(x-1)^{p^a},(x+v)^{p^a}\},\\ (x-1)^{p^a} & \mbox{ if } h(x)=(x-1)^{p^a}, \\ (x+u)^{p^a} & \mbox{ if } h(x)= (x+v)^{p^a}.\end{cases}$$
Then $\pi'$ is a desired map. $\square$.

\smallskip
{\bf Proof of Corollary \ref{cor1.4}.}  Let $\lambda=(\lambda_1,\cdots,\lambda_k)$ be the cycle type of $\sigma\in \Sigma_n$, and let $\sigma':=(1,2,\cdots,\lambda_1)\cdots$ $(\sum^{k-1}_{j=1}\lambda_j+1, \sum^{k-1}_{j=1}\lambda_j+2,\cdots, n)\in\Sigma_n$. Then $\sigma$ and $\sigma'$ are conjugate, and their permutation matrices are similar. Since similar matrices have the same minimal polynomial, we have $m_{c_\sigma}(x)=m_{c_{\sigma'}}(x)={\rm lcm}(x^{\lambda_1}-1,\cdots, x^{\lambda_k}-1)$, the least common multiple of all $x^{\lambda_i}-1$, $i\in [k]$.
According to Lemma \ref{per}, we have
$$\qquad (\dag) \quad \mathcal{E}_{c_{\sigma}} = \{ f(x)^{p^{\nu_p(\lambda_i)}}\mid i\in[k], f(x)~\mbox{is an irreducible factor of}~x^{\lambda_i}-1\}.$$
Note that $x-1\notin  \mathcal{E}_{c_{\sigma}}$ if and only if $\nu_p(\lambda_i)>0$ for all $i\in [k]$ if and only if $\sigma=s(\sigma)$ if and only $r(\sigma)=id$.

By definition, the cycle type of  $s(\sigma)$ is $(\lambda_{i_1},\cdots, \lambda_{i_t},\underbrace{1,\cdots, 1}_{k-t})$,  where $\{\lambda_{i_1},\cdots,\lambda_{i_t}\}=\{\lambda_i\mid i\in [k], \nu_p(\lambda_i)>0\}$ and the number of $1$-cycles of $s(\sigma)$ is exactly the number of points involved in $p$-regular cycles of $\sigma$. Thanks to $(\dag)$, we obtain
$$\mathcal{E}_{c_{s(\sigma)}}=\begin{cases}\{u(x)\in \mathcal{E}_{c_\sigma}\mid u(x)~\mbox{is reducible in}~ R[x]\} & \mbox{ if } s(\sigma)=\sigma,\\ \{u(x)\in \mathcal{E}_{c_\sigma}\mid u(x)~\mbox{is reducible in}~ R[x]\}\cup \{x-1\} & \mbox{ if } s(\sigma)\neq \sigma.\end{cases}$$
This implies $$\mathcal{R}_{c_{s(\sigma)}}=\begin{cases}\mathcal{R}_{c_\sigma}& \mbox{ if } s(\sigma)\neq id,\\ \emptyset & \mbox{ if } s(\sigma)= id.\end{cases}$$Thus, if $s(\sigma)\neq id$, then $P_{c_\sigma}((x-1)^{p^a})=P_{c_{s(\sigma)}}((x-1)^{p^a})$ and $P_{c_{s(\sigma)}}(h(x))=P_{c_\sigma}(h(x))\setminus\{1\}$ for $h(x)\in \mathcal{R}_{c_{s(\sigma)}}\setminus\{(x-1)^{p^a}\}$, where $(x-1)^{p^a}$ is the exceptional maximal elementary divisor of $c_\sigma$. Similarly, if $s(\tau)\neq id$, then $P_{c_\tau}((x-1)^{p^b})=P_{c_{s(\tau)}}((x-1)^{p^b})$ and $P_{c_{s(\tau)}}(g(x))=P_{c_\tau}(g(x))\setminus\{1\}$ for $g(x)\in \mathcal{R}_{c_{s(\tau)}}\setminus\{(x-1)^{p^b}\}$, where $(x-1)^{p^b}$ is the exceptional maximal elementary divisor of $c_\tau$.

By Theorem \ref{main1}(1), it suffices to show that if $c_{\sigma}\stackrel{S}\sim c_{\tau}$, then $c_{s(\sigma)}\stackrel{S}\sim c_{s(\tau)}$. Suppose $c_{\sigma}\stackrel{S}\sim c_{\tau}$. By Lemma \ref{iso-pr}(2), we have $s(\sigma)=id$ if and only if  $s(\tau)=id$, and so the statement is clear when $s(\sigma)=id$.

In the following, we assume that $s(\sigma)\neq id$ and $s(\tau)\neq id$, that is, both $\mathcal{R}_{c_\sigma}$ and $\mathcal{R}_{c_\tau}$ are nonempty. Then $p>0$  divides the orders of both $\sigma$ and $\tau$, and there is a bijection $\pi: \mathcal{R}_{c_\sigma}\to \mathcal{R}_{c_\tau}$ such that $R[x]/(h(x))\simeq R[x]/((h(x)\pi)$ as algebras and either $P_{c_\sigma}(h(x))=P_{c_\tau}((h(x))\pi)$ or $P_{c_\sigma}(h(x))=\mathcal{J}_{P_{c_\tau}((h(x))\pi)}$ for $h(x)\in \mathcal{R}_{c_\sigma}$.

According to $(\dag)$, $P_{c_\sigma}(h(x))$ and $P_{c_\tau}((h(x))\pi)$ consist of only $p$-powers. Thus, by Lemma \ref{p-p}, we obtain the equality $P_{c_\sigma}(h(x))=P_{c_\tau}((h(x))\pi)$ for all $h(x)\in \mathcal{R}_{c_\sigma}$.
By Lemma \ref{map}, we may assume $((x-1)^{p^a})\pi=(x-1)^{p^b}$. Then $$P_{c_{s(\sigma)}}((x-1)^{p^a})=P_{c_\sigma}((x-1)^{p^a})=P_{c_\tau}((x-1)^{p^b})=P_{c_{s(\tau)}}((x-1)^{p^b}).$$ Note that $P_{c_{s(\sigma)}}(h(x))=P_{c_\sigma}(h(x))\setminus\{1\}=P_{c_\tau}(h(x)\pi)\setminus\{1\}=P_{c_{s(\tau)}}(h(x)\pi)$ for $h(x)\in \mathcal{R}_{c_{s(\sigma)}}\setminus\{(x-1)^{p^a}\}$. Thus the restriction of $\pi$ to $\mathcal{R}_{c_{s(\sigma)}}$ gives rise to an $S$-equivalence between $s(\sigma)$ and $s(\tau)$. $\square$

\smallskip
The above proof and \cite[Theorem 1.1]{lx1} imply the strong conclusion.
\begin{Koro} Let $R$ be a field of characteristic $p>0$. If $\sigma\in \Sigma_n$ and $\tau\in \Sigma_m$ are $p$-singular permutations, then $S_n(c_{\sigma},R)$ and $S_m(c_{\tau},R)$ are stably equivalent if and only if they are Morita equivalent. \end{Koro}

In general, the converse of Corollary \ref{cor1.4} does not have to be true. This is illustrated by the example.

\begin{Bsp}\label{ex4.5} {\rm Let $R$ be an algebraically closed field of characteristic $3$. We take $\sigma\in \Sigma_8$ with the cycle type $(6,2)$, and $\tau\in \Sigma_7$ with the cycle type $(6,1)$. In this case, $s(\sigma)$ is a permutation of the cycle type $(6, 1^2)$, while $s(\tau)$ has the cycle type $(6,1)$. Clearly,  $S_8(c_{s(\sigma)},R)$ and $S_7(c_{s(\tau)},R)$ are stably equivalent by Corollary \ref{cor1.4}. In fact, they are Morita equivalent by \cite[Lemma 2.5(2)]{xz2}.

By Lemma \ref{per}, $\mathcal{E}_{c_\sigma}=\{ (x-1)^3, (x+1)^3, x-1,x+1 \}$, $\mathcal{R}_{c_\sigma}=\{(x-1)^3,(x+1)^3\}$, $P_{c_\sigma}((x-1)^3)=P_{c_\sigma}((x+1)^3)=\{1,3\}$; $\mathcal{E}_{c_\tau}=\{ (x-1)^3, (x+1)^3, x-1 \}$, $\mathcal{R}_{c_\tau}=\{(x-1)^3,(x+1)^3\}$, $P_{c_\tau}((x-1)^3)=\{1,3\}$ and $P_{c_\tau}((x+1)^3)=\{3\}$. Clearly, there are not any bijections between $\mathcal{R}_{c_\sigma}$ and $\mathcal{R}_{c_\sigma}$ such that $c_\sigma$ and $c_\tau$ are $S$-equivalent, and therefore $S_8(c_{\sigma},R)$ and $S_7(c_{\tau},R)$ cannot be stably equivalent by Theorem \ref{main1}(1).
}\end{Bsp}

\subsection{Homological dimensions and stable equivalences: Proof of Corollary \ref{cor1.5}\label{sect3.5}}
In this subsection, we show that stable equivalences preserve dominant, finitistic and global dimensions of centralizer matrix algebras over perfect fields. This proves Corollary \ref{cor1.5} .

\smallskip
We first point out a nice property of stable equivalences of centralizer matrix algebras over fields.

\begin{Prop}\label{S-SM} Let $R$ be a field, $c\in M_n(R)$ and $d\in M_m(R)$. Then $S_n(c,R)$ and $S_m(d,R)$ are stably equivalent if and only if their non-semisimple parts are stably equivalent of Morita type.
\end{Prop}

{\it Proof.} Suppose that $F$ is a stable equivalence between $S_n(c,R)$ and $S_m(d,R)$. Then $F$ preserves non-semisimple blocks by Lemma \ref{sta}. Thus, to complete the proof, it suffices to show that the corresponding non-semisimple blocks under $F$ are stably equivalent of Morita type. Let $A$ be an arbitrary non-semisimple block of $S_n(c,R)$, and assume that $B$ is the corresponding non-semsimple block of $S_m(d,R)$ via $F$. By our convention, $A=\End_{U_i}(M_i)$ for some $i\in [l_c]$ and $B=\End_{V_j}(N_j)$ for some $j\in [l_d]$. Then, by Lemma \ref{sta-down}, there is an isomorphism $U_i\simeq V_j$ of algebras. This enables us to view $N_j$ as an $U_i$-module. Recall that $M_i$ (respectively, $N_j$) is a faithful module which contains the left regular modules $U_i$ (respectively, $V_j$) as a direct summand. Then, applying Lemma \ref{Stc2}, we deduce that either $\mathcal{B}(M_i)_{\mathscr{P}}\simeq \mathcal{B}(N_j)_{\mathscr{P}}$ or $\mathcal{B}(M_i)_{\mathscr{P}}\simeq \Omega_{U_i}(\mathcal{B}(N_j)_{\mathscr{P}})$. Thus $\End_{U_i}(U_i\oplus \mathcal{B}(M_i)_{\mathscr{P}})$ and $\End_{V_j}(V_j\oplus \mathcal{B}(N_j)_{\mathscr{P}})$ are stably equivalent of Morita type in both cases (apply Lemma \ref{almm} for the second case). Hence $A=\End_{U_i}(M_i)$ and $B=\End_{V_j}(N_j)$ are stably equivalent of Morita type.  $\square$

\smallskip
Stable equivalences of Morita type between algebras over fields preserve the  global and finitistic dimensions of algebras (see Subsection \ref{sect2.3}). So we have the first two statements of the following result.

\begin{Koro} \label{homdim} Suppose that $R$ is a field, $c\in M_n(R)$ and $d\in M_m(R)$ such that $c$ and $d$ are $S$-equivalent. Then

$(1)$ $\gd\big(S_n(c,R)\big) = \gd\big(S_m(d,R)\big).$

$(2)$ $\fd\big(S_n(c,R)\big) = \fd\big(S_m(d,R)\big)$.

$(3)$ $\dm\big(S_n(c,R)\big) = \dm\big(S_m(d,R)\big)$.
\end{Koro}
{\it Proof.} Note that the global dimension and finitistic dimension of an algebra is the maximum of that of its non-semisimple blocks. Thus (1) and (2) follow immediately from Proposition \ref{S-SM}. It remains to prove (3). It was proved in \cite[Lemma 4.8(1)]{lx1} that the dominant dimension of a centralizer matrix algebra belongs to $\{2, \infty\}$. Thus it suffices to show that $\dm(S_n(c,R)) =\infty$ if and only if $\dm(S_m(d,R))=\infty$. By \cite[Lemma 4.8(2)]{lx1}, $\dm(S_n(c,R)) =\infty$ (respectively, $\dm(S_m(d,R)) =\infty$) if and only if all $P_c(f(x))$ (respectively, $P_d(g(x))$) for $f(x)\in \mathcal{M}_c$ (respectively, $g(x)\in \mathcal{M}_d$) are singleton sets. Note that the power indice set for any irreducible elementary divisor is $\{1\}$. Hence, by the definition of $S$-equivalences, we see that $P_c(f(x))$ for $f(x)\in \mathcal{M}_c$ are singleton sets if and only if  $P_d(g(x))$ for $g(x)\in \mathcal{M}_d$ are singleton sets. This completes the proof. $\square$

\begin{Rem}{\rm
(1) Under the assumptions that the ground field is algebraically closed and algebras considered have neither nodes nor semisimple direct summands, Corollary \ref{homdim} (1) and (3) were showed by Mart\'inez-Villa in \cite{MV2}. But Corollary \ref{homdim} removes all of these restrictions for centralizer matrix algebras.

(2) In general, stable equivalences do not have to preserve dominant, finitistic and global dimensions. This can be seen by the stable equivalence between $A:=\mathbb{Q}[x]/(x^2)$ and its node-eliminated algebra $A'$ which is isomorphic to the $2$ by $2$ upper triangular matrix algebra over $\mathbb{Q}$. Hence Corollary \ref{homdim} reflects special features of centralizer matrix algebras and also shows that the $2\times 2$ upper triangular matrix algebra over a field $R$ is the non-example of centralizer matrix algebra of the smallest dimension over $R$.
}
\end{Rem}

To end this section, we propose the following open problems.

\medskip
\textbf{Problem 1.} For a principal integral domain $R$, $c\in M_n(R)$ and $d\in M_m(R)$, we conjecture that $S_n(c,R)$ and $S_m(d,R)$ are stably equivalent if and only if $c$ and $d$ are $S$-equivalent.

\textbf{Problem 2.} What are invariants of $S$-equivalence of matrices over fields?

\textbf{Problem 3.} What are equivalent characterizations of $S$-equivalence in terms of matrix theory?

\medskip
\textbf{Acknowledgements.} The research work was partially supported by National Natural Science Foundation of China (Grants 12031014, 12350710787). Also, the corresponding author C. C. Xi thanks the
Tianyuan project of NSFC for partial support (Grant 12226314).

{\footnotesize
\smallskip
Xiaogang Li,
Shenzhen International Center for Mathematics, Southern University of Science and Technology, 518055
Shenzhen, Guangdong, P. R. China;

{\tt Email: 2200501002@cnu.edu.cn}

\smallskip
Changchang Xi,

School of Mathematical Sciences, Capital Normal University, 100048
Beijing, P. R. China, and

School of Mathematics and Statistics, Shaanxi Normal University, 710119 Xi'an, P. R. China

{\tt Email: xicc@cnu.edu.cn}
}

\end{document}